\documentclass{amsart}
\theoremstyle{plain}

\newtheorem{thm}{Theorem}[section]
\newtheorem{theorem}[thm]{Theorem}
\newtheorem{cor}[thm]{Corollary}

\newtheorem{lemma}[thm]{Lemma}

\newtheorem{prop}[thm]{Proposition}

\newtheorem*{conjecture*}{Conjecture}

\newtheorem*{question*}{Question}

\theoremstyle{definition}

\newtheorem*{definitions*}{Definitions}

\newtheorem*{rem*}{Remark}
\newtheorem{remark}[thm]{Remark}
\newtheorem*{remark*}{Remark}
\newtheorem{remarks}[thm]{Remarks}
\newtheorem*{remarks*}{Remarks}
\newtheorem*{example*}{Example}
\newtheorem{example}[thm]{Example}
\newtheorem*{examples*}{Examples}
\newtheorem{examples}[thm]{Examples}

\newtheorem*{notation*}{Notation}
\newtheorem*{convention*}{Convention}
\newtheorem*{conventions*}{Conventions}
\newtheorem*{note*}{Note}
\newtheorem*{exercise*}{Exercise}
\newtheorem*{bibliographical-note*}{Bibliographical note}

\newcommand{\zeroindent}{\parindent0cm \parskip1ex}
\newcommand{\acknowledgements}{{\em Acknowledgements. }}

\newenvironment{theoremlist}%
{\begin{list}{{\rm(\alph{enumi}) }}{\usecounter{enumi}

\leftmargin0cm \labelsep0cm \rightmargin0cm \topsep-1em
\setlength{\labelwidth}{\fill}}} {\end{list}}

{

\begin{enumerate}}{\end{enumerate}}

\newenvironment{alphalist}%
{

\begin{enumerate}}{\end{enumerate}}

\newenvironment{romanlist}%
{

\begin{enumerate} \parsep0cm
\itemsep0cm \parskip0cm}{\end{enumerate}}

{

\begin{enumerate}}{\end{enumerate}}

{

\begin{enumerate}}{\end{enumerate}}

\newenvironment{remarkslist}%
{\begin{list}{(\alph{enumi}) }{\usecounter{enumi}

\leftmargin0cm \labelsep0cm
\itemsep1em \rightmargin0cm
\setlength{\labelwidth}{\fill}}} {\end{list}}
\newcommand{\includefigure}[3]{%
\begin{figure}[#3]
\begin{center}
\epsfig{file=#2}\\ \caption{\label{fig:#1}}
\end{center}
\end{figure}}
\newcommand{\R}{\mathbb{R}}
\newcommand{\Z}{\mathbb{Z}}

\newcommand{\C}{\mathbb{C}}

\newcommand{\id}{\mathrm{id}}
\newcommand{\im}{\mathrm{im}}

\newcommand{\suchthat}{\; | \;}
\newcommand{\iso}{\cong}           
\newcommand{\htp}{\simeq}            
\newcommand{\smooth}{C^\infty}
\newcommand{\CP}[1]{\C {\mathrm P}^{#1}}
\newcommand{\RP}[1]{\R {\mathrm P}^{#1}}
\newcommand{\half}{{\textstyle\frac{1}{2}}}

\newcommand{\mo}{(M,\omega)}
\renewcommand{\o}{\omega}
\renewcommand{\O}{\Omega}
\newcommand{\Aut}{\mathrm{Aut}}

\newcommand{\Symp}{\mathrm{Sp}}
\newcommand{\leftsc}{\langle}
\newcommand{\rightsc}{\rangle}
\newcommand{\Diff}{\mathrm{Diff}}

\usepackage{xspace}
\usepackage[matrix,arrow]{xy}
\usepackage[dvips]{epsfig}
\CompileMatrices

\zeroindent \numberwithin{equation}{section} \theoremstyle{definition}

\theoremstyle{definition}
\newtheorem{examples-remarks}[thm]{Examples and comments}

\renewcommand{\paragraph}[1]{{\bf #1.}}

\newcommand{\tSymp}{\widetilde{\Symp}}

\newcommand{\LL}{\mathcal{L}}
\newcommand{\tPhi}{\tilde{\Phi}}

\renewcommand{\det}{\mathrm{det}}
\newcommand{\tLambda}{\tilde{\Lambda}}
\newcommand{\tlambda}{\tilde{\lambda}}
\newcommand{\tphi}{\tilde{\phi}}
\newcommand{\tL}{\tilde{L}}
\newcommand{\Lag}{\mathrm{Lag}}
\newcommand{\Laggr}{\Lag^{\mathrm{gr}}}
\newcommand{\Autgr}{\Aut^{\mathrm{gr}}}

\newcommand{\B}{\mathcal{B}}
\newcommand{\J}{\mathbf{J}}
\newcommand{\gen}[1]{\leftsc #1 \rightsc}
\newcommand{\moduli}{{\mathcal M}}
\newcommand{\cotangent}{T^*}
\newcommand{\ucotangent}{T_1^*}
\newcommand{\var}{\mathrm{var}}
\newcommand{\tchi}{\tilde{\chi}}
\newcommand{\tpsi}{\tilde{\psi}}

\newcommand{\HP}[1]{\mathbb{H} \mathrm{P}^{#1}}
\newcommand{\tI}{\tilde{I}}
\newcommand{\re}{\mathrm{re}\;}
\newcommand{\ttau}{\tilde{\tau}}
\newcommand{\isotopic}{\htp}
\newcommand{\grisotopic}{\isotopic_{\mathrm{gr}}}
\newcommand{\JJ}{\mathcal{J}}
\newcommand{\tH}{\tilde{H}}
\newcommand{\tSigma}{\tilde{\Sigma}}
\newcommand{\diag}{\mathrm{diag}}

\title{Graded Lagrangian submanifolds}
\author{Paul Seidel}
\date{Revised version, July 5, 1999.}

\begin{document}
\maketitle

\section{Introduction}

Floer theory assigns, in favourable circumstances, an abelian group $HF(L_0,L_1)$ to a
pair $(L_0,L_1)$ of Lagrangian submanifolds of a symplectic manifold $\mo$. This group
is a qualitative invariant, which remains unchanged under suitable deformations of
$L_0$ or $L_1$. Following Floer \cite{floer88e} one can equip $HF(L_0,L_1)$ with a
canonical relative $\Z/N$-grading, where $1 \leq N \leq \infty$ is a number which
depends on $\mo$, $L_0$ and $L_1$ (for $N = \infty$ we set $\Z/N = \Z$). {\em Relative}
mostly means that the grading is unique up to an overall shift, although there are also
cases with more complicated behaviour. In this paper we take a different approach to
the grading: we consider Lagrangian submanifolds equipped with certain extra structure
(these are what we call graded Lagrangian submanifolds). This extra structure removes
the ambiguity and defines an {\em absolute} $\Z/N$-grading on Floer cohomology. There
is also a parallel notion of graded symplectic automorphism, which bears the same
relation to the corresponding version of Floer theory. Both concepts were first
discovered by Kontsevich, at least for $N = \infty$; see \cite[p.\ 134]{kontsevich94}.
Somewhat later, the present author came upon them independently.

One way to approach the definition of graded Lagrangian submanifold is to start with
the case $N = 2$. It is well-known that orientations of $L_0$ and $L_1$ determine an
absolute $\Z/2$-grading $HF(L_0,L_1) = HF^0(L_0,L_1) \oplus HF^1(L_0,L_1)$. One can
reformulate this as follows: consider the natural fibre bundles $\LL,\LL^{or}
\longrightarrow M$ whose fibres are the unoriented resp. oriented Lagrangian
Grassmannians of the tangent spaces $TM_x$. Any Lagrangian submanifold $L$ comes with a
canonical section $s_L: L \longrightarrow \LL|L$, and an orientation of $L$ is the same
as a lift of this section to $\LL^{or}$. Hence the right objects for a Floer theory
with an absolute $\Z/2$-grading are pairs $(L,\tL)$ consisting of a Lagrangian
submanifold and a lift $\tL: L \longrightarrow \LL^{or}$ of $s_L$. In order to define
the absolute $\Z/N$-grading one proceeds in the same way, only that $\LL^{or}$ must be
replaced by a $\Z/N$-covering of $\LL$ of a certain kind. Such coverings, which we call
Maslov coverings, need not exist in general, and they are also not unique. In fact,
choosing an $N$-fold Maslov covering is equivalent to lifting the structure group of
$TM$ from $\Symp(2n)$ to a certain finite extension $\Symp^N(2n)$; and the particularly
simple situation for $N = 2$ is due to the fact that $\Symp^2(2n) \iso \Symp(2n) \times
\Z/2$.

In itself this `graded symplectic geometry' is not particularly deep, but it does make
Floer cohomology into a more powerful invariant. To put it bluntly, the advantage of
the new framework is this: in passing to graded Lagrangian submanifolds there is a
choice of $\Z/N$ for any Lagrangian submanifold $L$ (the choice of the lift of $s_L$).
In comparison, if one uses only the relative grading, there is a $\Z/N$-ambiguity for
any {\em pair} of Lagrangian submanifolds, and this greater amount of choice entails a
loss of information. We illustrate this through three applications, which form the main
part of this paper.

{\bf (a) Lagrangian submanifolds of $\CP{n}$.} We prove that any Lagrangian submanifold
$L \subset \CP{n}$ must satisfy $H^1(L;\Z/(2n+2)) \neq 0$ (the actual result is
slightly sharper, see Theorem \ref{th:embedding}).

{\bf (b) Symplectically knotted Lagrangian spheres.} The paper \cite{seidel98b}
provides examples of compact symplectic four-manifolds (with boundary) $M$ with the
following property: there is a family of embedded Lagrangian two-spheres $L^{(k)}
\subset M$, $k \in \Z$, such that any two of them are isotopic as smooth submanifolds,
but no two are isotopic as Lagrangian submanifolds. In such a situation we say that $M$
contains infinitely many symplectically knotted Lagrangian two-spheres. The examples in
\cite{seidel98b} were constructed using a special class of symplectic automorphisms,
called generalized Dehn twists, and the main step in the proof was a Floer cohomology
computation using Pozniak's \cite{pozniak} Morse-Bott type spectral sequence. Both the
construction and the proof can be generalized to produce Lagrangian $n$-spheres with
the same property for all even $n$.

Here, using the method of graded Lagrangian submanifolds, we will first reprove the
result from \cite{seidel98b} and its generalization in a considerably simpler way.
Then, by a more complicated construction, we produce similar examples of Lagrangian
$n$-spheres for all odd $n \geq 5$. The reason why the remaining case $n = 3$ cannot be
settled in the same way is topological, and seems to have nothing to do with Floer
theory.

We can also improve on \cite{seidel98b} in a different direction, by showing that
suitable $K3$ and Enriques surfaces contain infinitely many symplectically knotted
Lagrangian two-spheres. These are the first known examples of closed symplectic
manifolds with this property. As a by-product one obtains that for these manifolds the
map $\pi_0(\Aut\mo) \longrightarrow \pi_0(\Diff(M))$ has infinite kernel, sharpening a
result of \cite{seidel97}. Unfortunately, at the present state of development in Floer
theory, it is impossible for technical reasons to carry out a similar argument in
dimensions $>4$.

{\bf (c) Weighted homogeneous singularities.} Let $p \in \C[x_0,\dots,x_n]$, $n \geq
1$, be a weighted homogeneous polynomial with an isolated critical point at the origin.
One can introduce the Milnor fibre of $p$, which is a compact symplectic manifold
$(M^{2n},\o)$ with boundary, and the symplectic monodromy $f \in \Aut(M,\partial M,\o)$
of the Milnor fibration. This refines the usual notion of geometric monodromy by taking
into account the symplectic geometry of the situation. We will show that $[f] \in
\pi_0(\Aut(M,\partial M,\o))$ has infinite order whenever the sum of the weights is
$\neq 1$. It is not known whether the condition on the weights is really necessary.

It should be mentioned (although this will not be used later on) that this application
and the previous one are related. In fact, generalized Dehn twists are maps modelled on
the monodromy of the quadratic singularity $p(x) = x_0^2 + \dots + x_n^2$, and the
construction of odd-dimensional knotted Lagrangian spheres is inspired by the monodromy
of the singularity $p(x) = x_0^2 + \dots + x_{n-1}^2 + x_n^3$ of type $(A_2)$.

The importance of `graded symplectic geometry' for these applications varies. For (a)
and (c) its role is that of a convenient language. In fact one could replace it by
monodromy considerations in the style of \cite{seidel96} without changing the essence
of the argument. For (a) there is also a more algebraic argument, based on the fact
that $HF(L,L)$ is a module over the quantum cohomology $QH^*(\CP{n})$. The situation in
(b) is different, since the `graded' framework allows us to state a basic geometric
property of generalized Dehn twists (Lemma \ref{th:dehn-shift}) which it seems hard to
encode in any other way.

{\em Notation and conventions.} All manifolds are usually assumed to be connected. The
automorphism group of a symplectic manifold $\mo$ will be denoted by $\Aut\mo$. If $M$
is compact, we equip this group with the $\smooth$-topology. If $M$ is a symplectic
manifold with nonempty boundary, $\Aut(M,\partial M,\o) \subset \Aut\mo$ is the
subgroup of automorphisms $\phi$ which are equal to the identity on some neighbourhood
(depending on $\phi$) of the boundary. Lagrangian submanifolds are always assumed to be
compact; if the symplectic manifold has a boundary, any Lagrangian submanifold is
assumed to lie in the interior. $\Lag\mo$ stands for the space of Lagrangian
submanifolds of $M$, with the $\smooth$-topology. $S^1$ will often be identified with
$\R/\Z$. (Co)homology groups have $\Z$-coefficients unless otherwise stated.

\acknowledgements This paper is an offshoot of my joint work with Mikhail Khovanov;
several of the ideas presented here arose in conversations with him. I am indebted to
Maxim Kontsevich for explaining his joint work with Fukaya to me. The idea of
considering symplectic manifolds with a circle action on the boundary, which is the
subject of section \ref{sec:class}, was suggested to me by Leonid Polterovich and Yakov
Eliashberg. This paper was written while I was staying at the Max Planck Institute
(Bonn) and ETH Z{\"u}rich, and I would like to thank both institutions for their
hospitality. I am indebted to the referee for several helpful suggestions.

\section{Basic notions}
\subsection{Linear algebra\label{subsec:linear-algebra}}

By a $\Z/N$-covering ($1 \leq N \leq \infty$) of a space $X$ we mean a covering $X^N$
with covering group $\Z/N$. Such coverings are classified up to isomorphism by
$H^1(X;\Z/N)$. For connected $X$, this correspondence associates to a homomorphism
$\pi: \pi_1(X) \longrightarrow \Z/N$ the covering $X^N = \tilde{X} \times_\pi \Z/N$,
where $\tilde{X}$ is the universal cover. If $X$ is a connected Lie group, all
$\Z/N$-coverings of it (even the non-connected ones) have canonical Lie group
structures.

Let $(V^{2n},\beta)$ be a symplectic vector space, $\Symp(V,\beta)$ the linear
symplectic group, and $\LL(V,\beta)$ the Lagrangian Grassmannian, which parametrizes
linear Lagrangian subspaces of $V$. Both $\Symp(V,\beta)$ and $\LL(V,\beta)$ are
connected with infinite cyclic fundamental group. Moreover, there are preferred
generators $\delta(V,\beta) \in H^1(\Symp(V,\beta))$ and $C(V,\beta) \in
H^1(\LL(V,\beta))$ (the second one is called the Maslov class) so that one can
canonically identify the fundamental groups with $\Z$. $\Symp(V,\beta)$ acts
transitively on $\LL(V,\beta)$, and any orbit is a map $\Symp(V,\beta) \longrightarrow
\LL(V,\beta)$ which takes $C(V,\beta)$ to $2\delta(V,\beta)$. For $1 \leq N \leq
\infty$, let $\LL^N(V,\beta)$ be the $\Z/N$-covering of $\LL(V,\beta)$ which
corresponds to the image of $C(V,\beta)$ in $H^1(\LL(V,\beta);\Z/N)$. The $\Z/N$-action
on $\LL^N(V,\beta)$ will be denoted by $\rho$. Define $\Symp^N(V,\beta)$ to be the
group of pairs $(\Phi,\tPhi)$ consisting of $\Phi \in \Symp(V,\beta)$ and a
$\Z/N$-equivariant diffeomorphism $\tPhi$ of $\LL^N(V,\beta)$, which is a lift of the
action of $\Phi$ on $\LL(V,\beta)$. This is a Lie group and fits into an exact sequence
\[
1 \longrightarrow \Z/N \longrightarrow \Symp^N(V,\beta) \longrightarrow \Symp(V,\beta)
\longrightarrow 1,
\]
where $\Z/N$ is the central subgroup of pairs $(\Phi,\tPhi) = (\id,\rho(k))$. It
follows that $\Symp^N(V,\beta)$ must be isomorphic to some $\Z/N$-covering of
$\Symp(V,\beta)$. The next Lemma identifies that covering.

\begin{lemma} \label{th:sympn}
$\Symp^N(V,\beta)$ is isomorphic (as a Lie group) to the $\Z/N$-covering of
$\Symp(V,\beta)$ associated to the image of $2\delta(V,\beta)$ in
$H^1(\Symp(V,\beta);\Z/N)$.
\end{lemma}

\proof Let $\tSymp(V,\beta)$ be the universal cover of $\Symp(V,\beta)$. Take a loop
$\phi: [0;1] \longrightarrow \Symp(V,\beta)$ with $\phi(0) = \phi(1) = \id$ and
$\leftsc \delta(V,\beta), [\phi] \rightsc = 1$, and let $\tphi(1) \in \tSymp(V,\beta)$
the endpoint of the lift $\tphi$ of $\phi$ with $\tphi(0) = \id$. The action of
$\Symp(V,\beta)$ on $\LL(V,\beta)$ can be lifted uniquely to an action of
$\tSymp(V,\beta)$ on $\LL^N(V,\beta)$. This action commutes with the $\Z/N$-action
$\rho$, and $\tphi(1)$ acts in the same way as $\rho(2)$. Therefore one obtains a
homomorphism
\[
\tSymp(V,\beta) \times_{\pi} \Z/N \longrightarrow \Symp^N(V,\beta),
\]
where $\pi: \pi_1(\Symp(V,\beta)) = \Z \longrightarrow \Z/N$ is multiplication by two.
It is not difficult to see that this is an isomorphism, which proves the desired
result. \qed

As an example consider the case $N = 2$. One can identify $\LL^2(V,\beta)$ with the
oriented Lagrangian Grassmannian $\LL^{or}(V,\beta)$. Since $\Symp(V,\beta)$ acts
naturally on $\LL^{or}(V,\beta)$ one has $\Symp^2(V,\beta) \iso \Z/2 \times
\Symp(V,\beta)$. More generally, one can try to compare $\Symp^N(V,\beta)$ with the
more obvious covering $\Symp^N(V,\beta)'$ of $\Symp(V,\beta)$ obtained from the mod $N$
reduction of $\delta(V,\beta)$. One finds that $\Symp^N(V,\beta)' \iso
\Symp^N(V,\beta)$ if $N$ is finite and odd, and that
\begin{equation} \label{eq:comparison}
\Symp^{2N}(V,\beta) \iso \Symp^N(V,\beta)' \times_{\Z/N} \Z/2N.
\end{equation}
In particular, $\Symp^N(V,\beta)$ is connected iff $N$ is finite and odd, and has two
connected components otherwise.

Let $J$ be a $\beta$-compatible complex structure on $V$, and $g$ the corresponding
inner product. Recall that the unitary group $U(V,J,g) \subset \Symp(V,\beta)$ is a
deformation retract, and that $\delta(V,\beta)$ is represented by the determinant
$U(V,J,g) \longrightarrow S^1$. Let $U^N(V,J,g)$ be the $\Z/N$-covering of $U(V,J,g)$
determined by the mod $N$ reduction of $2\delta(V,\beta)$. These coverings are clearly
deformation retracts of $\Symp^N(V,\beta)$, and they are explicitly given by
\begin{equation} \label{eq:un}
U^N(V,J,g) = \left\{\;\;\begin{matrix} \{(\Phi,q) \in U(V,J,g) \times S^1 \suchthat
\det(\Phi)^2 = q^N\} & \text{for } N <\infty,\\ \{(\Phi,t) \in U(V,J,g) \times \R
\suchthat \det(\Phi)^2 = e^{2\pi i t}\} & \text{for } N = \infty.\end{matrix}\right.
\end{equation}

In future we will abbreviate $\LL(\R^{2n},\o_{\R^{2n}})$ by $\LL(2n)$. Similarly we
will write $C(2n)$, $\LL^N(2n)$, $\Symp(2n)$, $\Symp^N(2n)$, and $U^N(n)$.

\subsection{Graded symplectic geometry\label{subsec:graded}}

Let $(M^{2n},\o)$ be a symplectic manifold, possibly with boundary. Let $P
\longrightarrow M$ be the principal $\Symp(2n)$-bundle associated to the symplectic
vector bundle $(TM,\o)$, and $\LL \longrightarrow M$ the natural fibre bundle whose
fibres are the Lagrangian Grassmannians $\LL_x = \LL(TM_x,\o_x)$. One can identify $\LL
= P \times_{\Symp(2n)} \LL(2n)$.

\begin{alphalist}
\item \label{item:spn-structures}
An $\Symp^N(2n)$-structure ($1 \leq N \leq \infty$) on $M$ is a principal
$\Symp^N(2n)$-bundle $P^N \longrightarrow M$ together with an isomorphism $P^N
\times_{\Symp^N(2n)} \Symp(2n) \iso P$.

\item
An {\em $N$-fold Maslov covering} is a $\Z/N$-covering $\LL^N \longrightarrow \LL$
whose restriction to $\LL_x = \LL(TM_x,\o_x)$, for any $x \in M$, is isomorphic to
$\LL^N(TM_x,\o_x)$. The $\Z/N$-action on $\LL^N$ will always be denoted by $\rho$.

\item
A global Maslov class mod $N$ is a class $C^N \in H^1(\LL;\Z/N)$ whose restriction to
any fibre is the mod $N$ reduction of the ordinary Maslov class.
\end{alphalist}

There are canonical bijections between (isomorphism classes of) these three kinds of
objects. If $P^N$ is an $\Symp^N(2n)$-structure then the associated fibre bundle with
fibre $\LL^N(2n)$ is an $N$-fold Maslov covering. Conversely, in the presence of a
Maslov covering $\LL^N$, the transition maps of any system of local trivializations of
$(TM,\o)$ have canonical lifts from $\Symp(2n)$ to $\Symp^N(2n)$ which satisfy the
cocycle condition, hence define an $\Symp^N(2n)$-structure. The connection between
Maslov coverings and global Maslov classes is obvious. Now assume that we have chosen
an $\o$-compatible almost complex structure $J$ on $M$, and consider the line bundle
$\Delta = \Lambda^n(TM,J)^{\otimes 2}$.

\begin{alphalist}
\setcounter{enumi}{3}
\item \label{item:determinant}
An $N$-th root ($1 \leq N <\infty$) of $\Delta$ is a complex line bundle $Z$ together
with an isomorphism $r: Z^{\otimes N} \longrightarrow \Delta$. Two pairs $(Z,r)$,
$(Z',r')$ are called equivalent if there is an isomorphism $j: Z \longrightarrow Z'$
such that $r'j = r$. In addition, we define an $\infty$-th root to be a trivialization
of $\Delta$, and two of them are called equivalent iff they are homotopic.
\end{alphalist}

There is a canonical bijection between $\Symp^N(2n)$-structures and equivalence classes
of $N$-th roots of $\Delta$; it is defined as follows. Let $P_U$ be the principal
$U(n)$-bundle associated to $(TM,\o,J)$. A $U^N(n)$-structure on $M$ is a principal
$U^N(n)$-bundle together with an isomorphism of the associated $U(n)$-bundle with
$P_U$. Because $U^N(n) \subset \Symp^N(2n)$ is a deformation retract, there is a
canonical bijection between $U^N(n)$-structures and $\Symp^N(2n)$-structures. On the
other hand, by looking at \eqref{eq:un} one sees that a $U^N(n)$-structure is just a
choice of $N$-th root of $\Delta$. Among the equivalent notions
\ref{item:spn-structures}--\ref{item:determinant} we will most frequently work with
Maslov coverings, since that is convenient for dealing with Lagrangian submanifolds.

\begin{lemma} \label{th:newbasic-one}
$(M^{2n},\o)$ admits an $N$-fold Maslov covering iff $2c_1\mo$ goes to zero in
$H^2(M;\Z/N)$. The isomorphism classes of such coverings (provided that any exist) form
an affine space over $H^1(M;\Z/N)$.
\end{lemma}

\proof This is immediate if one uses an almost complex structure and the description
\ref{item:determinant}. Alternatively one can use \ref{item:spn-structures} and an
argument based on the exact sequence
\[
0 \rightarrow H^1(M;\Z/N) \rightarrow H^1(M;\Symp^N(2n)) \rightarrow H^1(M;\Symp(2n))
\rightarrow H^2(M;\Z/N)
\]
of non-abelian cohomology groups, just as in the classification of spin structures in
\cite[Appendix A]{lawson-michelsohn}. \qed

Let $\LL^N$ be an $N$-fold Maslov covering on $\mo$. For every Lagrangian submanifold
$L \subset M$ there is a natural section $s_L: L \longrightarrow \LL|L$, $s_L(x) = TL_x
\in \LL(TM_x,\o_x)$. An {\em $\LL^N$-grading} of $L$ is a lift $\tL: L \longrightarrow
\LL^N$ of $s_L$. The pair $(L,\tL)$ is called an {\em $\LL^N$-graded Lagrangian
submanifold}. We write $\Laggr(M,\o;\LL^N)$ for the set of such pairs, and equip it
with the topology which comes from the space of compact submanifolds of $\LL^N$ (by
considering the image of $\tL$). This topology defines the notion of an isotopy of
graded Lagrangian submanifolds. Clearly, if $\tL$ is an $\LL^N$-grading of $L$ then so
is $\rho(k) \circ \tL$ for any $k \in \Z/N$. This defines a free $\Z/N$-action on
$\Laggr(M,\o;\LL^N)$.

\begin{lemma} \label{th:newbasic-two}
The forgetful map $\Laggr(M,\o;\LL^N) \longrightarrow \Lag(M,\o)$ is a $\Z/N$-covering
of its image. The image itself consists of those $L$ such that $s_L^*(C^N) \in
H^1(L;\Z/N)$ is zero, where $C^N$ is the global Maslov class corresponding to $\LL^N$.
In particular, a Lagrangian submanifold with $H^1(L;\Z/N) = 0$ always admits an
$\LL^N$-grading.
\end{lemma}

\proof By definition, $L$ admits an $\LL^N$-grading iff $s_L^*(\LL^N) \longrightarrow
L$ is a trivial covering, which is equivalent to $s_L^*(C^N) = 0$. The rest is obvious.
\qed

Let $\LL^N$ be an $N$-fold Maslov covering on $\mo$ and $\phi$ a symplectic
automorphism. There is a natural map $\phi^\LL: \LL \longrightarrow \LL$,
$\phi^\LL(\Lambda) = D\phi(\Lambda)$, which covers $\phi$. An {\em $\LL^N$-grading} of
$\phi$ is a $\Z/N$-equivariant diffeomorphism $\tphi$ of $\LL^N$ which is a lift of
$\phi^\LL$. The pair $(\phi,\tphi)$ is called an {\em $\LL^N$-graded symplectic
automorphism}. Such pairs form a group which we denote by $\Autgr(M,\o;\LL^N)$. If $M$
is compact, we equip this group with the topology induced from embedding it into
$\Diff(\LL^N)^{\Z/N}$. The pairs $(\phi,\tphi) = (\id,\rho(k))$ form a central subgroup
$\Z/N \subset \Autgr(M,\o;\LL^N)$. $\LL^N$-graded symplectic automorphisms act
naturally on $\LL^N$-graded Lagrangian submanifolds by $(\phi,\tphi)(L,\tL) =
(\phi(L),\tphi \circ \tL \circ \phi^{-1})$.

\begin{lemma} \label{th:newbasic-three}
Let $\LL^N$ be an $N$-fold Maslov covering. The forgetful homomorphism
$\Autgr(M,\o;\LL^N) \longrightarrow  \Aut\mo$ fits into an exact sequence
\[
1 \longrightarrow \Z/N \longrightarrow \Autgr(M,\o;\LL^N) \longrightarrow \Aut\mo
\stackrel{\partial}{\longrightarrow} H^1(M;\Z/N).
\]
Here $\partial$ is not a group homomorphism, but it satisfies $\partial(\phi\psi) =
\psi^*\partial(\phi) + \partial(\psi)$, so that its kernel is a subgroup of $\Aut\mo$.
If $M$ is compact then this is a sequence of topological groups, with $\Z/N$ and
$H^1(M;\Z/N)$ discrete.
\end{lemma}

\proof By definition, a symplectic automorphism $\phi$ admits an $\LL^N$-grading iff
the two Maslov coverings $(\phi^\LL)^*(\LL^N)$ and $\LL^N$ are isomorphic. By Lemma
\ref{th:newbasic-one} the difference between these two coverings can be measured by a
class in $H^1(M;\Z/N)$. We define $\partial(\phi)$ to be this class. The rest is easy.
\qed

\begin{remark} \label{th:boundary}
Assume that $M$ has nonempty boundary. Then one can consider the subgroup
$\Autgr(M,\partial M,\o;\LL^N) \subset \Autgr(M,\o;\LL^N)$ consisting of pairs with
$\phi \in \Aut(M,\partial M,\o)$ and where $\tphi|\LL^N_x = \id$ for any $x \in
\partial M$ (here $\LL^N_x$ is the part of $\LL^N$ which lies over $\LL_x$). The
central elements $(\id,\rho(k))$, $k \neq 0$, do not lie in $\Autgr(M,\partial
M,\o;\LL^N)$. In fact one has an exact sequence (with $\partial$ defined in a similar
way as before)
\[
1 \longrightarrow \Autgr(M,\partial M,\o;\LL^N) \longrightarrow \Aut(M,\partial M,\o)
\stackrel{\partial}{\longrightarrow} H^1(M,\partial M;\Z/N).
\]
\end{remark}

The minimal Chern number $N_M$ of $\mo$ is defined to be the positive generator of the
group $\leftsc c_1(M), H_2(M) \rightsc \subset \Z$. Similarly, the relative minimal
Chern number $N_L$ of a Lagrangian submanifold $L \subset M$ is the positive generator
of the group $\leftsc 2c_1(M,L), H_2(M,L) \rightsc$, where now $2c_1(M,L) \in H^2(M,L)$
is the relative first Chern class. These numbers, or variants of them, are familiar
from the definition of the relative grading on Floer cohomology. Their relationship to
the concepts introduced here, at least in the case $H_1(M) = 0$, is as follows.

\begin{lemma} \label{th:newbasic-four}
A symplectic manifold $\mo$ with $H_1(M) = 0$ admits a Maslov covering $\LL^N$ of order
$N$ iff $N$ divides $2N_M$. Moreover, this covering is unique up to isomorphism. A
Lagrangian submanifold $L \subset M$ admits an $\LL^N$-grading iff $N$ divides $N_L$.
\end{lemma}

\proof Because $H_1(M) = 0$, it follows from the universal coefficient sequence that
$2c_1(M)$ goes to zero in $H^2(M;\Z/N)$ iff $\leftsc 2c_1(M),x \rightsc$ is a multiple
of $N$ for any $x \in H_2(M)$. In view of Lemma \ref{th:newbasic-one}, this proves the
first part.

Take a compact oriented surface $\Sigma$ and a map $w: (\Sigma, \partial \Sigma)
\longrightarrow (M,L)$. The number $\leftsc 2c_1(M,L), [w] \rightsc \in \Z$ can be
computed as follows: choose a trivialization of the pullback $w^*\LL$, that is to say,
a fibre bundle map
\[
\xymatrix{{\LL(2n) \times \Sigma} \ar[rr]^\tau \ar[d] && {\LL} \ar[d] \\ {\Sigma}
\ar[rr]^w && {M.}}
\]
One has $\tau^{-1} \circ s_L \circ (w|\partial\Sigma)(x) = (\lambda(x),x)$ for some map
$\lambda:
\partial\Sigma \longrightarrow \LL(2n)$, and $\leftsc 2c_1(M,L),[w] \rightsc = \leftsc C(2n),
\lambda_*[\partial\Sigma]\rightsc$.

Let $\LL^N$ be the unique Maslov covering of some order $N$ on $\mo$ and $C^N$ its
global Maslov class. The pullback $\tau^*(C^N) \in H^1(\LL(2n) \times \Sigma;\Z/N)$ is
of the form $C(2n) + y$ for some $y \in H^1(\Sigma;\Z/N)$. Hence in $\Z/N$ one has
\begin{align*}
& \leftsc s_L^*(C^N), w_*[\partial\Sigma] \rightsc = \\ &= \leftsc C^N, (s_L \circ
w|\partial\Sigma)_*[\partial\Sigma] \rightsc = \leftsc \tau^*(C^N), (\tau^{-1} \circ
s_L \circ w|\partial\Sigma )_*[\partial\Sigma] \rightsc =\\ &= \leftsc C(2n),
\lambda_*[\partial\Sigma] \rightsc + \leftsc y, [\partial\Sigma] \rightsc = \leftsc
2c_1(M,L), [w] \rightsc.
\end{align*}
If $N|N_L$ then $\leftsc 2c_1(M,L), [w] \rightsc$ is always a multiple of $N$. Since
one can choose $w$ in such a way that $w_*[\partial\Sigma]$ is an arbitrary class of
$H_1(L)$, it follows that $s_L^*(C^N) = 0$, which means that $L$ admits an
$\LL^N$-grading. The converse is equally simple. \qed

In future we will use the following notation. Instead of $(\phi,\tphi)$ and $(L,\tL)$
we will often write only $\tphi$ and $\tL$. The action of $\Autgr(M,\o;\LL^N)$ on
$\Laggr(M,\o;\LL^N)$ will be written as $\tphi(\tL)$. We will denote $(\id,\rho(-k))
\in \Autgr(M,\o;\LL^N)$ by $[k]$ and call it the $k$-fold shift operator. The graded
Lagrangian submanifold $\rho(-k) \circ \tL$, which is obtained from $\tL$ by the action
of $[k]$, will be denoted by $\tL[k]$. The similarity with homological algebra is
intentional, and the sign in the definition of $[k]$ has been introduced with that in
mind.

\subsection{Examples\label{subsec:examples}}

We will now complement the basic definitions by several examples and remarks, some of
which will be used later on.

\begin{example} \label{ex:orientation}
Since $\Symp^2(2n) \iso \Symp(2n) \times \Z/2$, an $\Symp^2(2n)$-structure is just the
choice of a real line bundle $\xi$ on $M$. The corresponding two-fold Maslov covering,
which we denote by $\LL^{or,\xi}$, is the space of pairs $(\Lambda,o)$, where $\Lambda
\in \LL$ is a Lagrangian subspace of $TM_x$ and $o$ is an orientation of the vector
space $\Lambda \otimes_\R \xi_x$. An $\LL^{or,\xi}$-grading of a Lagrangian submanifold
$L \subset M$ is the same as an orientation of $TL \otimes (\xi|L)$. An
$\LL^{or,\xi}$-grading of a symplectic automorphism $\phi$ is the same as a bundle
isomorphism $\phi^*(\xi) \longrightarrow \xi$. In particular, the trivial line bundle
yields a two-fold Maslov covering $\LL^{or}$ for which a grading of a Lagrangian
submanifold is just an orientation, and such that $\Autgr(M,\o;\LL^{or}) \iso \Aut\mo
\times \Z/2$.
\end{example}

\begin{example} \label{ex:spin}
One can associate to any spin structure on $M$ an $\Symp^4(2n)$-structure. The reason,
in the notation of \eqref{eq:comparison}, is that the restriction of the universal
cover of $GL^+(2n)$ to the subgroup $\Symp(2n)$ is again a nontrivial double cover,
hence isomorphic to $\Symp^2(2n)'$, and that $\Symp^4(2n) \iso \Symp^2(2n)'
\times_{\Z/2} \Z/4$. Note that not all $\Symp^4(2n)$-structures arise in this way.
\end{example}

\begin{example} \label{ex:calabi-yau}
The following discussion relates our point of view to an another one (which is also
originally due to Kontsevich). Let $(V,\beta)$ be a $2n$-dimensional symplectic vector
space, $J$ a compatible complex structure, and $g$ the corresponding inner product. Set
$\Delta(V,J) = \Lambda^n(V,J)^{\otimes 2}$, and let $S\Delta(V,J) \subset \Delta(V,J)$
be the unit circle (with respect to the metric induced by $g$). One can define a
fibration with simply-connected fibres
\[
\det^2: \LL(V,\beta) \longrightarrow S\Delta(V,J)
\]
by $\det^2(\Lambda) = (e_1 \wedge \dots \wedge e_n)^{\otimes 2}$, where $(e_j)$ is any
orthonormal basis of $(\Lambda,g|\Lambda)$. After choosing an element $\Theta \in
\Delta(V,J)^*$ of unit length one can identify $S\Delta(V,J)$ with $S^1$. In this way
one obtains a map $\det^2_\Theta: \LL(V,\beta) \longrightarrow S^1$. The Maslov class
$C(V,\beta)$ is equal to the pullback of the standard generator $[S^1]$. Hence
$\LL^\infty(V,\beta)$ is isomorphic to the pullback of the standard covering $\R
\longrightarrow S^1$.

Now let $(M,\o)$ be a symplectic manifold and $J$ a compatible almost complex
structure. Assume that $2c_1\mo = 0$, which means that $\Delta =
\Lambda^n(TM,J)^{\otimes 2}$ is trivial. Choose a section $\Theta$ of $\Delta^*$ (in
other words, a quadratic complex $n$-form) which has length one everywhere. As before
this determines a map $\det_\Theta^2: \LL \longrightarrow S^1$, and one can define an
$\infty$-fold Maslov covering by
\begin{equation} \label{eq:theta-covering}
\LL^\infty = \{ (\Lambda,t) \in \LL \times \R \suchthat \det^2_\Theta(\Lambda) =
e^{2\pi i t}\}.
\end{equation}
An $\LL^\infty$-grading of a Lagrangian submanifold $L$ is just a lift of the map
$\det_\Theta^2 \circ s_L : L \longrightarrow S^1$ to $\R$. This approach is
particularly useful in complex geometry: let $(M,\o,J)$ be a Calabi-Yau manifold, take
a covariantly constant holomorphic $n$-form $\theta$ of unit length, and set $\Theta =
\theta^{\otimes 2}$. A Lagrangian submanifold $L \subset M$ is called special if
$(\im\; \theta)|L = 0$. This condition is equivalent to $\det_\Theta^2 \circ s_L \equiv
1 \in S^1$. It follows that special Lagrangian submanifolds have a canonical
$\LL^\infty$-grading.
\end{example}

\begin{example} \label{ex:lagrangian-distribution}
Let $\mo$ be a symplectic manifold which admits a Lagrangian distribution $S \subset
TM$. Then one can define an $\infty$-fold Maslov covering $\LL^\infty$ simply by
putting together the universal covers of $\LL(TM_x,\o_x)$ with base point $S_x$, for
all $x$. The global Maslov class of this covering is represented by the fibrewise
Maslov cycle
$
I = \bigcup_{x \in M} \{ \Lambda \in \LL_x \suchthat \Lambda \cap S_x \neq 0\}
$
with its canonical co-orientation. Any Lagrangian submanifold which is either tangent
or transverse to $S$ admits an $\LL^\infty$-grading (there is even a preferred one).
Typical examples are cotangent bundles $M = \cotangent X$ with either the vertical
distribution (tangent spaces along the fibres) or the horizontal one (with respect to
some Riemannian metric on $X$).
\end{example}

\begin{example} \label{ex:surfaces}
Let $M$ be an oriented closed surface of genus $g \geq 2$, and $\o$ a volume form on
it. Using Lemma \ref{th:newbasic-one} one can see that $\mo$ admits many different
Maslov coverings of any order which divides $4g-4$. However, in a sense there is no
really good choice of covering:

\begin{prop}
For every Maslov covering $\LL^N$ of order $N > 2$ on $M$, there is an automorphism
$\phi \in \Aut\mo$ which does not admit an $\LL^N$-grading.
\end{prop}

\includefigure{surface}{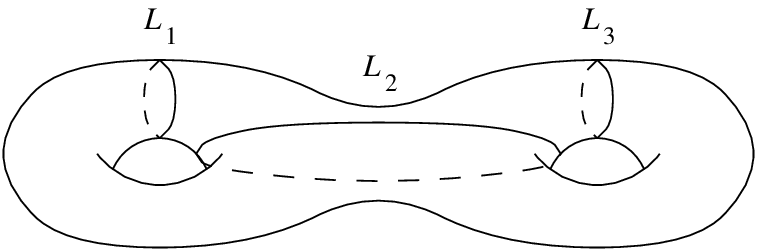}{ht}

\proof Choose a Maslov covering $\LL^N$ and let $C^N$ be the corresponding global
Maslov class. One can associate to any oriented embedded curve $L \subset M$ a number
$R(\LL^N,L) = \leftsc s_L^*(C^N), [L] \rightsc \in \Z/N$. $L$ admits an $\LL^N$-grading
iff it has zero rotation number, and a symplectic automorphism $\phi \in \Aut\mo$
admits an $\LL^N$-grading iff $R(\LL^N,L) = R(\LL^N,\phi(L))$ for all curves $L$. We
need two more facts, whose proofs we leave to the reader: (a) if $\Sigma \subset M$ is
a surface whose boundary is formed by the curves $L_1,\dots,L_k$ then $\sum_j
R(\LL^N,L_j) \equiv 2\chi(\Sigma)$ mod $N$. (b) If $t_L$ is the positive Dehn twist
along a curve $L$, then $R(\LL^N,t_L(L')) = R(\LL^N,L') - ([L'] \cdot [L]) R(\LL^N,L)$
for any $L'$. Now consider the curves $L_k$ shown in fig. \ref{fig:surface}. By (a) we
have $R(\LL^N,L_1) + R(\LL^N,L_2) + R(\LL^N,L_3) = -2$. Hence there is a $\nu \in
\{1,2,3\}$ such that $R(\LL^N,L_\nu) \neq 0$. Take a curve $L'$ with $[L'] \cdot
[L_\nu] = \pm 1$. Property (b) shows that $R(\LL^N,t_{L_\nu}(L')) \neq R(\LL^N,L')$,
which means that $t_{L_\nu}$ does not admit an $\LL^N$-grading. This is for $g = 2$,
but one sees immediately that the argument generalizes to all $g > 2$. \qed

In contrast, for $M = T^2$ there is exactly one $\infty$-fold Maslov covering
$\LL^\infty$ (namely, the one coming from the standard trivialization of $TM$) with the
property that any $\phi \in \Aut\mo$ admits an $\LL^\infty$-grading.
\end{example}

\subsection{Two kinds of index\label{subsec:index}}

Let $(V,\beta)$ be a $2n$-dimensional symplectic vector space. The {\em Maslov index
for paths} \cite{viterbo87} \cite{robbin-salamon93} assigns a half-integer
$\mu(\lambda_0,\lambda_1) \in \half\Z$ to any pair of paths $\lambda_0,\lambda_1: [a;b]
\longrightarrow \LL(V,\beta)$. $\mu(\lambda_0,\lambda_1)$ is an integer iff
$\dim(\lambda_0(a) \cap \lambda_1(a)) \equiv \dim(\lambda_0(b) \cap \lambda_1(b))$ mod
$2$. We will now adapt this invariant to our situation. Fix some $1 \leq N \leq
\infty$. Let $\tLambda_0,\tLambda_1 \in \LL^N(V,\beta)$ be a pair of points whose
images in $\LL(V,\beta)$ intersect transversely. Choose two paths
$\tlambda_0,\tlambda_1 : [0;1] \longrightarrow \LL^N(V,\beta)$ with $\tlambda_0(0) =
\tlambda_1(0)$ and $\tlambda_0(1) = \tLambda_0$, $\tlambda_1(1) = \tLambda_1$. Let
$\lambda_0, \lambda_1$ be the projections of these paths to $\LL(V,\beta)$. The {\em
absolute Maslov index} of $(\tLambda_0,\tLambda_1)$ is defined by
\[
\tilde{\mu}(\tLambda_0,\tLambda_1) = n/2 - \mu(\lambda_0,\lambda_1)  \in \Z/N.
\]
One can easily show that this is independent of all choices. From the standard
properties of $\mu$ \cite{robbin-salamon93} one derives the following facts:

\begin{romanlist}
\item \label{item:difference} Let $\lambda_0,\lambda_1: [a;b] \longrightarrow
\LL(V,\beta)$ be two paths such that $\lambda_0(a) \cap \lambda_1(a) = \lambda_0(b)
\cap \lambda_1(b) = 0$. Lift them to paths $\tlambda_0,\tlambda_1$ in $\LL^N(V,\beta)$.
Then
\[
\tilde{\mu}(\tlambda_0(b),\tlambda_1(b)) - \tilde{\mu}(\tlambda_0(a),\tlambda_1(a)) =
-\mu(\lambda_0,\lambda_1).
\]
\item \label{item:shifting} $\tilde{\mu}(\rho(k)\tLambda_0,\rho(l)\tLambda_1) =
\tilde{\mu}(\tLambda_0,\tLambda_1) - k + l$.

\item \label{item:automorphisms}
For any $(\Phi,\tilde{\Phi}) \in \Symp^N(V,\beta)$ one has
$\tilde{\mu}(\tilde{\Phi}(\tLambda_0), \tilde{\Phi}(\tLambda_1)) =
\tilde{\mu}(\tLambda_0,\tLambda_1)$.

\item \label{item:symmetry} $\tilde{\mu}(\tLambda_1,\tLambda_0) = n -
\tilde{\mu}(\tLambda_0,\tLambda_1)$.

\item \label{item:morse} Let $\Lambda,\Lambda^\perp \subset V$ be two complementary
Lagrangian subspaces. Let $B$ be a nondegenerate quadratic form on $\Lambda$, and $A:
\Lambda \longrightarrow \Lambda^\perp$ the unique linear map such that $B(v) =
\beta(v,Av)$. Take the path $\lambda: [0;1] \longrightarrow \LL(V,\beta)$ given by
$\lambda(t) =  \mathrm{graph}(tA)$, and lift it to a path $\tlambda$ in
$\LL^N(V,\beta)$. Then $\tilde{\mu}(\tlambda(1),\tlambda(0))$ is the Morse index of $B$
(mod $N$).

\item For $N = 2$, if one identifies $\LL^2(V,\beta) = \LL^{or}(V,\beta)$ then $(-1)^{\tilde{\mu}}$
agrees with the intersection number of oriented Lagrangian subspaces up to a constant
$(-1)^{n(n+1)/2}$.
\end{romanlist}

The application to graded symplectic geometry is as follows. Let $\mo$ be a symplectic
manifold with an $N$-fold Maslov covering $\LL^N$, and $(L_0,\tL_0)$, $(L_1,\tL_1)$ a
pair of $\LL^N$-graded Lagrangian submanifolds which intersect transversally. Then one
can associate to any point $x \in L_0 \cap L_1$ an absolute index mod $N$,
\[
\tilde{I}(\tL_0,\tL_1;x) \stackrel{\mathrm{def}}{=} \tilde{\mu}(\tL_0(x),\tL_1(x)).
\]

The {\em Conley-Zehnder index} associates to any path $\phi: [a;b] \longrightarrow
\Symp(V,\beta)$ a half-integer $\zeta(\phi) \in \half\Z$. It can be reduced to the
Maslov index for paths as follows: take $(V',\beta') = (V,-\beta) \oplus (V,\beta)$ and
consider the two paths in $\LL(V',\beta')$ given by $\lambda_0(t) =
\mathrm{graph}(\phi(t))$ and $\lambda_1(t) = \Delta$ (the diagonal). Then $\zeta(\phi)
= \mu(\lambda_0,\lambda_1)$ (here we are following \cite[Remark 5.4]{robbin-salamon93};
it seems that the definition in \cite{salamon-zehnder92} has the opposite sign). Now
let $(\Phi,\tilde{\Phi}) \in \Symp^N(V,\beta)$ be a point such that $\det(1-\Phi) \neq
0$. Choose a path $\tilde{\phi}: [0;1] \longrightarrow \Symp^N(V,\beta)$ from
$(\id,\rho(k)) \in \Symp^N(V,\beta)$, for some $k \in \Z/N$, to $(\Phi,\tilde{\Phi})$,
and project it to a path $\phi$ in $\Symp(V,\beta)$. The {\em absolute Conley-Zehnder
index} is $\tilde{\zeta}(\Phi,\tilde{\Phi}) = n - \zeta(\phi) - k \in \Z/N$. This is
independent of all choices and has the following properties:

\begin{romanlist}
\item \label{item:cz-conjugation} it is invariant under conjugation.

\item \label{item:cz-shifting} $\tilde{\zeta}(\Phi, \tilde{\Phi} \circ \rho(l)) =
\tilde{\zeta}(\Phi,\tilde{\Phi}) - l$.

\item \label{item:cz-inverse}
$\tilde{\zeta}(\Phi^{-1},\tilde{\Phi}^{-1}) = 2n - \tilde{\zeta}(\Phi,\tilde{\Phi})$.

\item \label{item:cz-morse}
Let $B$ be a nondegenerate quadratic form on $V$, and $A: V \longrightarrow V$ the
linear map given by $\o(x,Ax) = B(x)$. Take the path $\phi(t) = \exp(tA)$ in
$\Symp(V,\beta)$ and lift it to a path $\tilde{\phi}$ in $\Symp^N(V,\beta)$ with
$\tilde{\phi}(0) = (\id,\id)$. Then $\tilde{\zeta}(\phi(t),\tilde{\phi}(t))$ is equal
to the Morse index of $B$ (mod $N$) for sufficiently small $t>0$.
\end{romanlist}

Given a symplectic manifold $\mo$, an $N$-fold Maslov covering $\LL^N$, and an
$\LL^N$-graded symplectic automorphism $(\psi,\tpsi)$ which has nondegenerate fixed
points, one can associate to any fixed point $x$ an absolute index $\tilde{Z}(\tpsi;x)
= \tilde{\zeta}(D\psi_x,\tpsi(x)) \in \Z/N$.

\subsection{Lagrangian surgery\label{subsec:surgery}}

We will now discuss an example which shows that the absolute Maslov index appears even
in elementary questions about graded Lagrangian submanifolds.

Take an embedding $\gamma: \R \longrightarrow \C$ with $\gamma(t) = t$ for $t \leq
-1/2$, $\gamma(t) = it$ for $t \geq 1/2$, and $\gamma(\R) \cap -\gamma(\R) =
\emptyset$. Then $H = \bigcup_{t \in \R} \gamma(t)S^{n-1} \subset \C^n$ is a Lagrangian
submanifold with respect to the standard symplectic form. Outside the unit ball $B^{2n}
\subset \C^n$ one has $H \cap (\C^n \setminus B^{2n}) = (\R^n \cup i\R^n) \setminus
B^{2n}$. $H$ is called a {\em Lagrangian handle}. It is used in the following way: let
$(M^{2n},\o)$ be a symplectic manifold and $L_1,L_2 \subset M$ a pair of Lagrangian
submanifolds which intersect transversely in a single point $\{x_0\} = L_1 \cap L_2$.
There is always an embedding $j: B^{2n} \longrightarrow M$ with $j(0) = x_0$,
$j^{-1}(L_1) = \R^n \cap B^{2n}$, $j^{-1}(L_2) = i\R^n \cap B^{2n}$, and $j^*\o =
\epsilon\,\o_{\C^n}$ for some $\epsilon>0$. Then one can form the embedded connected
sum
\[
L_1 \# L_2 = (L_1 \setminus \im(j)) \cup (L_2 \setminus \im(j)) \cup j(H \cap B^{2n})
\]
which is again a Lagrangian submanifold. This process, which is independent of all
choices up to Lagrangian isotopy, is usually called {\em Lagrangian surgery}. It has
been studied by Polterovich \cite{polterovich91} and others. Our conventions are those
of \cite[Appendix]{seidel98b} and differ from Polterovich's. Note that Lagrangian
surgery is not symmetric: $L_1 \# L_2$ and $L_2 \# L_1$ are the two possibilities of
resolving the self-intersection of $L_1 \cup L_2$ (this can be seen clearly already in
the case $n = 1$).

Take the $\infty$-fold Maslov covering $\LL^\infty$ on $\C^n$ induced by the quadratic
complex $n$-form $\Theta = (dz_1 \wedge \dots \wedge dz_n)^{\otimes 2}$ as in
\eqref{eq:theta-covering}. An $\LL^\infty$-grading of a Lagrangian submanifold $L
\subset \C^n$ is the same as a map $\tL: L \longrightarrow \R$ such that $e^{2\pi i
\tL(x)} = \det^2_\Theta(TL_x)$ for all $x$. For $L_1 = \R^n$ and $L_2 = i\R^n$ we
choose the gradings $\tL_1 \equiv 0$, $\tL_2 \equiv 1 - n/2$. Then the absolute index
at the origin is $\tilde{I}(\tL_1,\tL_2;0) = 1$.

\begin{lemma} There is an $\LL^\infty$-grading $\tH$ of the Lagrangian handle $H$ which agrees with
$\tL_1$ on $\R^n \setminus B^{2n}$ and with $\tL_2$ on $i\R^n \setminus B^{2n}$.
\end{lemma}

\proof Take $y_0 \in S^{n-1}$, and let $y_0^\perp \subset \R^n$ be its orthogonal
complement. The tangent space of $H$ at a point $y = \gamma(t)y_0$ is $TH_y = \R
\gamma'(t) y_0 \oplus \gamma(t) y_0^\perp$. Therefore
\begin{equation} \label{eq:determinant}
\det^2_\Theta(TH_y) =
\frac{\gamma'(t)^2\gamma(t)^{2n-2}}{|\gamma'(t)^2\gamma(t)^{2n-2}|} \in S^1.
\end{equation}
As $t$ goes from $-\infty$ to $\infty$, $\gamma(t)^2/|\gamma(t)|^2$ makes half a turn
clockwise from $1$ to $-1$, and $\gamma'(t)^2/|\gamma'(t)|^2$ makes half a turn
counterclockwise from $1$ to $-1$. It follows that one can find a map $\alpha \in
\smooth(\R,\R)$ with $\alpha(t) = 0$ for $t \leq -1/2$ and $\alpha(t) = 1 - n/2$ for $t
\geq 1/2$, such that $e^{2\pi i \alpha(t)}$ equals the r.h.s. of
\eqref{eq:determinant}. Then $\tilde{H}(e^{it}y_0) = \alpha(t)$ is a grading of $H$
with the desired property. \qed

From these local considerations one immediately obtains the following graded version of
Lagrangian surgery.

\begin{lemma} \label{th:graded-surgery}
Let $\mo$ be a symplectic manifold with an $\infty$-fold Maslov covering $\LL^\infty$.
Let $(L_1,\tL_1)$ and $(L_2,\tL_2)$ be two $\LL^\infty$-graded Lagrangian submanifolds
which intersect transversely and in a single point $x_0 \in M$. If
$\tilde{I}(\tL_1,\tL_2;x_0) = 1$ then the surgery $\Sigma = L_1 \# L_2$ has an
$\LL^\infty$-grading $\tilde{\Sigma}$ which agrees with $\tL_1$ on $\Sigma \cap L_1$,
and with $\tL_2$ on $\Sigma \cap L_2$. \qed
\end{lemma}

\subsection{Floer cohomology\label{subsec:floer}}

We can now introduce the absolute grading on Floer cohomology. The exposition in this
section is formal, in the sense that the conditions which are necessary to make Floer
cohomology well-defined will be suppressed. Our justification is that there is no
relation between these conditions and the problem of grading. Concretely, this means
that if the ordinary Floer cohomology $HF(L_0,L_1)$ of two Lagrangian submanifolds is
defined, and $L_0,L_1$ admit gradings $\tL_0,\tL_1$, then the graded version
$HF^*(\tL_0,\tL_1)$ is also defined. The discussion of the properties of
$HF^*(\tL_0,\tL_1)$ should be understood in the same way: they hold in the same
generality as their ungraded analogues.

Let $L_0,L_1$ be a pair of transversely intersecting Lagrangian submanifolds in a
symplectic manifold $(M^{2n},\o)$, and $\J = (J_t)_{0 \leq t \leq 1}$ a smooth family
of $\o$-compatible almost complex structures. For any two points $x_-,x_+ \in L_0 \cap
L_1$, let $\B(x_-,x_+)$ be the set of smooth maps $u: \R \times [0;1] \longrightarrow
M$ with $u(\R \times \{j\}) \subset L_j$ for $j = 0,1$, and $\lim_{s \rightarrow
\pm\infty} u(s,\cdot) = x_\pm$. Consider the subspace $\moduli(x_-,x_+;\J) \subset
\B(x_-,x_+)$ of maps which satisfy Floer's equation $\partial_su + J_t(u)\partial_tu =
0$. It has a natural action of $\R$ by translation in the first variable. For generic
choice of $\J$, $\moduli(x_-,x_+;\J)$ has a natural structure of a smooth manifold,
with connected components of different dimensions. Let $\moduli_k(x_-,x_+;\J)$ be the
$k$-dimension part of $\moduli(x_-,x_+;\J)$. In the simplest situation, such as in the
original work of Floer and in that of Oh \cite{oh93}, suitable assumptions on $\mo$ and
$L_0,L_1$ ensure that the quotients $\moduli_1(x_-,x_+;\J)/\R$ are finite sets. Then,
writing $n(x_-,x_+;\J) \in \Z/2$ for the number of points mod $2$ in
$\moduli_1(x_-,x_+;\J)/\R$, one defines a chain group $(CF(L_0,L_1),\partial_\J)$ as
follows: $CF(L_0,L_1)$ is the $\Z/2$-vector space freely generated by the points of
$L_0 \cap L_1$, and $\partial_\J\gen{x_+} = \sum_{x_-} n(x_-,x_+;\J) \gen{x_-}$. One
finds that $\partial_\J^2 = 0$, and the Floer cohomology is $HF(L_0,L_1;\J) = \ker
\partial_\J/ \im\;\partial_\J$ (this is Floer {\em co}homology because we have
exchanged the usual roles of $x_-$ and $x_+$ in the definition of $\partial_\J$).

Let $\LL^N$ be an $N$-fold Maslov covering on $M$, and assume that $L_0,L_1$ admit
$\LL^N$-gradings $\tL_0,\tL_1$. For $i \in \Z/N$, let $CF^i(\tL_0,\tL_1) \subset
CF(L_0,L_1)$ be the subspace generated by elements $\gen{x}$ where
$\tilde{I}(\tL_0,\tL_1;x) = i$. This defines a $\Z/N$-grading on $CF(L_0,L_1)$. Floer's
index theorem \cite{floer88e} together with property \ref{item:difference} of the
absolute Maslov index implies that $\partial_\J$ has degree one. Hence there is an
induced grading on Floer cohomology. We will refer to the $\Z/N$-graded group
$HF^*(\tL_0,\tL_1;\J)$ as {\em graded Floer cohomology}.

Let $\mathbf{L} = (L_t)_{0 \leq t \leq 1}$ be an exact isotopy of Lagrangian
submanifolds (exact means that it can be embedded in a Hamiltonian isotopy of $M$) and
$L$ a Lagrangian submanifold which intersects both $L_0$ and $L_1$ transversely. For
any almost complex structures $\J^-$, $\J^+$ such that $HF(L,L_0;\J^-)$ and
$HF(L,L_1;\J^+)$ are defined, there is a canonical isomorphism
\begin{equation} \label{eq:continuation}
q(\mathbf{L},\J^-,\J^+): HF(L,L_0;\J^-) \longrightarrow HF(L,L_1;\J^+).
\end{equation}
Now assume that we have a Maslov covering $\LL^N$, and that $L_0$ and $L$ admit
$\LL^N$-gradings $\tL_0,\tL$. Then the isotopy $(L_t)$ can be lifted to an isotopy
$(\tL_t)$ of $\LL^N$-graded Lagrangian submanifolds, and the map
\eqref{eq:continuation} has degree zero with respect to the gradings of
$HF^*(\tL,\tL_0;\J^-)$ and $HF^*(\tL,\tL_1;\J^+)$. To complete the construction of
graded Floer cohomology one follows the usual strategy: first, using the isomorphisms
\eqref{eq:continuation} for constant isotopies, one shows that graded Floer cohomology
is independent of the choice of almost complex structure (we will therefore omit $\J$
from the notation from now on). Secondly, using again \eqref{eq:continuation} but this
time for $C^1$-small isotopies, one extends the definition of graded Floer cohomology
to Lagrangian submanifolds which do not intersect transversally. Clearly, this extended
definition is still invariant under exact isotopies of graded Lagrangian submanifolds
(in both variables). Some other properties of graded Floer cohomology are: the shifting
formula
\begin{align*}
HF^j(\tL_0[k],\tL_1[l]) &\iso HF^{j-k+l}(\tL_0,\tL_1),\\ \intertext{invariance under
graded symplectic automorphisms,} HF^j(\tphi(\tL_0),\tphi(\tL_1)) &\iso
HF^j(\tL_0,\tL_1),\\ \intertext{and Poincar{\'e} duality:} HF^j(\tL_1,\tL_0) &\iso
HF^{n-j}(\tL_0,\tL_1)^{\vee}.
\end{align*}
These follow immediately from the properties of $\tilde{\mu}$ listed in section
\ref{subsec:index}. We also need to mention the isomorphism between $HF(L,L)$ and the
ordinary cohomology $H^*(L;\Z/2)$. This isomorphism holds in Floer's original setup
\cite{floer89}, but it can fail in more general situations, see \cite{oh96}. Whenever
it holds, so does the graded version
\[
HF^j(\tL,\tL) \iso \bigoplus_{i \in \Z} H^{j+iN}(L;\Z/2).
\]
The other version of Floer cohomology, that for symplectic automorphisms, will not be
used in this paper. Nevertheless, it seems appropriate to outline briefly the parallel
story about the grading. Let $\phi$ be an automorphism of $\mo$ which has
non-degenerate fixed points. One considers a chain group $CF(\phi)$ which has one basis
element $\gen{x}$ for any fixed point $x \in M$ of $\phi$, and a boundary operator
$\partial_\J$ defined as before, but this time using the space of maps $u: \R^2
\longrightarrow M$ which satisfy
\[
\begin{cases} & \partial_su + J_t(u)\partial_tu = 0, \\
& u(s,t) = \phi(u(s,t+1)), \end{cases}
\]
where $\J = (J_t)_{t \in \R}$ is a family of compatible almost complex structures
satisfying a suitable periodicity condition. If $\mo$ has a Maslov covering $\LL^N$ and
$\phi$ an $\LL^N$-grading $\tphi$ then, using the absolute Conley-Zehnder index
$\tilde{\zeta}$, one can define a $\Z/N$-grading on $CF(\phi)$. The index theorem of
\cite[section 4]{salamon-zehnder92} implies that $\partial_\J$ has degree one.
Therefore one obtains a $\Z/N$-graded Floer cohomology group $HF^*(\tphi)$. These
groups are invariant under Hamiltonian isotopies and satisfy $HF^*(\tphi \circ [k]) =
HF^{*-k}(\tphi)$, $HF^*(\tphi^{-1}) = HF^{2n-*}(\phi)^{\vee}$, $HF^*(\id) = H^*(M)$.
For the special case of monotone symplectic manifolds, a more extensive account of this
kind of Floer theory (without mention of the absolute grading) can be found in
\cite{dostoglou-salamon94}.

\section{Lagrangian submanifolds of $\CP{n}$}

In this section we use graded Floer cohomology to obtain some restrictions on the
topology of Lagrangian submanifolds of $\CP{n}$. Note that, by starting with the
familiar embeddings of $\RP{n}$ and $T^n$ and applying Lagrangian surgery, one can
construct many different Lagrangian submanifolds of $\CP{n}$.

\begin{theorem} \label{th:embedding}
Any Lagrangian submanifold $L \subset \CP{n}$ satisfies \begin{theoremlist}
\item $H^1(L;\Z/(2n+2)) \neq 0$,
\item $H^1(L;\Z/(2n+2)) \not\iso (\Z/2)^g$ for any $g \geq 2$,
\item if $H^1(L;\Z/(2n+2)) \iso \Z/2$ then $H^i(L;\Z/2) \iso \Z/2$ for all $i = 0, \dots, n$.
\end{theoremlist}
\end{theorem}

By Lemma \ref{th:newbasic-four}, $\CP{n}$ admits a unique Maslov covering $\LL^N$ of
any order $N$ which divides $2N_M = 2n+2$. Consider the Hamiltonian circle action given
by $\sigma(t) = \diag(e^{2\pi i t},1,\dots,1) \in U(n+1)$. One can lift $\sigma$
uniquely to a map $\tilde{\sigma}: [0;1] \longrightarrow \Autgr(M,\o;\LL^N)$ with
$\tilde{\sigma}(0) = \id$. By looking at any fixed point one sees that
$\tilde{\sigma}(1) = [-2]$. It follows that every $\LL^N$-graded Lagrangian submanifold
$\tL$ is graded Lagrangian isotopic to $\tL[-2]$, by an isotopy which is also exact.
This implies that $HF^*(\tL,\tL)$, whenever defined, is periodic with period two. From
this fact we will derive Theorem \ref{th:embedding}.

Before explaining the proof in detail, we need to recall the Floer cohomology for monotone
Lagrangian submanifolds as developed by Oh. The basic references are \cite{oh93} and
\cite{oh96}. Technical issues are discussed further in \cite{oh96c} and \cite{oh-kwon96}; see
also \cite{lazzarini98}. We will present the theory in a slightly simplified form. Let
$(M^{2n},\o)$ be a closed symplectic manifold which is monotone, that is, $[\o] = \lambda\,
c_1\mo$ for some $\lambda>0$. For any Lagrangian submanifold $L \subset M$, let $N_L \geq 1$
be the number defined in section \ref{subsec:graded}. Oh shows that $HF(L_0,L_1)$ is
well-defined, and invariant under Lagrangian isotopy, for all pairs $(L_0,L_1)$ such that
$N_{L_0}, N_{L_1} \geq 3$ and $H^1(L_0;\R) = H^1(L_1;\R) = 0$. Moreover one has

\begin{theorem}[Oh \protect{\cite[Theorem 5.1]{oh96}}] \label{th:oh}
Let $L \subset M$ be a Lagrangian submanifold with $H^1(L;\R) = 0$ and $N_L \geq 3$.
\begin{theoremlist} \item \label{item:oh-one}
If $N_L \geq n+2$ then $HF(L,L) \iso \bigoplus_i H^i(L;\Z/2)$.
\item \label{item:oh-two}
If $N_L = n+1$ then $HF(L,L)$ is either $\bigoplus_i H^i(L;\Z/2)$ or $\bigoplus_{i \neq
0,n} H^i(L;\Z/2)$.
\end{theoremlist}
\end{theorem}

The proof of this goes as follow \cite[p. 332]{oh96}. Let $H$ be a Morse function on
$L$, and $L'$ a small Lagrangian perturbation of $L$ constructed using $H$ and a
Darboux chart. We may assume that $H$ has only one local minimum $x_-$ and local
maximum $x_+$. The intersection points of $L$ and $L'$ are the critical points of $H$.
One can write the boundary operator on $CF(L,L')$ as $\partial_{\J} =
\partial_0 + \partial_1 + \dots$, where $\partial_k$ takes  critical points of Morse index $i$ to those
of Morse index $i+1-kN_L$. Floer \cite{floer89} proved that, for a suitable choice of
$\J$, $\partial_0$ can be identified with the boundary operator in a Morse cohomology
complex for $H$. Therefore the homology of $(CF(L,L'),\partial_0)$, which is sometimes
called the local Floer cohomology of $L$, is always isomorphic to $H^*(L;\Z/2)$. If
$N_L \geq n+2$ then for dimension reasons $\partial_k = 0$ for all $k>0$, which proves
\ref{item:oh-one}. If $N_L = n+1$ then $\partial_k = 0$ for $k \geq 2$,
$\partial_1\gen{x} = 0$ for all $x \neq x_+$, and $\partial_1\gen{x_+}$ can be either
zero or $\gen{x_-}$. This leads to the two possibilities in \ref{item:oh-two}.

Now let $\LL^N$ be a Maslov covering of order $N \geq 3$ on $M$. Lemma
\ref{th:newbasic-four} shows that any $\LL^N$-graded Lagrangian submanifold $(L,\tL)$
automatically satisfies $N_L \geq N \geq 3$. Hence the graded Floer cohomology groups
$HF^*(\tL_0,\tL_1)$ are well-defined for $\LL^N$-graded Lagrangian submanifolds with
zero first Betti number. Moreover, as a look at the proof shows, the obvious graded
analogue of Theorem \ref{th:oh} holds.

\proof[Proof of Theorem \ref{th:embedding}] (a) Assume that $L \subset \CP{n}$ ($n \geq
2$) is a Lagrangian submanifold with $H^1(L;\Z/(2n+2)) = 0$. This implies that
$H^1(L;\R) = 0$. By Lemma \ref{th:newbasic-two}, $L$ admits a grading $\tL$ with
respect to the unique Maslov covering $\LL^{2n+2}$ of order $2n+2$. Hence the graded
Floer cohomology $HF^*(\tL,\tL)$ is well-defined. Lemma \ref{th:newbasic-four} shows
that $N_L \geq 2n+2$, and by applying Theorem \ref{th:oh}\ref{item:oh-one} one finds
that
\begin{equation} \label{eq:the-gap}
HF^i(\tL,\tL) = \begin{cases} H^i(L;\Z/2) & 0 \leq i \leq n, \\ 0 & n+1 \leq i \leq
2n+1.
\end{cases}
\end{equation}
As discussed above, the circle action $\sigma$ on $\CP{n}$ provides a graded Lagrangian
isotopy between $\tL$ and $\tL[-2]$, which implies that $HF^*(\tL,\tL) \iso
HF^*(\tL,\tL[-2]) = HF^{*-2}(\tL,\tL)$. This is a contradiction, since
\eqref{eq:the-gap} is not two-periodic.

(b,c) Let $L \subset \CP{n}$ $(n \geq 2)$ be a Lagrangian submanifold with
$H^1(L;\Z/(2n+2)) = (\Z/2)^g$ for some $g \geq 1$. This implies that $H^1(L;\R) = 0$,
that $H^1(L;\Z/2) \iso (\Z/2)^g$, and that the homomorphism $\rho: H^1(L;\Z/(2n+2))
\longrightarrow H^1(L;\Z/(n+1))$ induced by the projection $\Z/(2n+2) \longrightarrow
\Z/(n+1)$ is zero.

Let $C^{2n+2} \in H^1(\LL;\Z/(2n+2))$ and $C^{n+1} \in H^1(\LL;\Z/(n+1))$ be the global
Maslov classes of the Maslov coverings $\LL^{2n+2},\LL^{n+1}$ on $\CP{n}$. Clearly
$C^{n+1}$ is obtained from $C^{2n+2}$ by reducing mod $(n+1)$. We conclude that
$s_L^*(C^{n+1}) = \rho(s_L^*(C^{2n+2})) = 0$, which means that $L$ admits an
$\LL^{n+1}$-grading $\tL$. The same argument as before shows that the $\Z/(n+1)$-graded
Floer cohomology $HF^*(\tL,\tL)$ is two-periodic. Lemma \ref{th:newbasic-four} says
that $N_L \geq n+1$. By Theorem \ref{th:oh}\ref{item:oh-two} there are two
possibilities for the Floer cohomology. One is that $HF^*(\tL,\tL) = H^*(L;\Z/2)$. In
this case it follows that $H^*(L;\Z/2)$, with the grading reduced mod $n+1$, is
two-periodic. Since $H^0(L;\Z/2) = \Z/2$ and $H^n(L;\Z/2) = \Z/2$, the periodicity
leads to $H^i(L;\Z/2) = \Z/2$ for all $i$, which proves both (b) and (c). It remains to
consider the other possibility, which is
\[
HF^i(\tL,\tL) = \begin{cases} H^i(L;\Z/2) & 0 < i < n, \\ 0 & i = 0,n. \end{cases}
\]
This contradicts the two-periodicity, because $HF^1(\tL,\tL) = H^1(L;\Z/2) = (\Z/2)^g$
while $HF^{-1}(\tL,\tL) = 0$. Hence this possibility cannot occur. \qed

\section{A class of symplectic automorphisms\label{sec:class}}

\subsection{The basic result\label{subsec:class}}

Let $(M,\o,\alpha)$ be a compact symplectic manifold with contact type boundary, and
$(\phi_t^K)$ an $S^1$-action on $\partial M$ which preserves $\alpha$. This means that
the symplectic form $\o \in \Omega^2(M)$ and the contact one-form $\alpha \in
\Omega^1(\partial M)$ are related by $d\alpha = \o|\partial M$, and that the Reeb
vector field $R$ of $\alpha$ satisfies $\o(N,R)>0$, where $N$ is any vector field
pointing outwards along $\partial M$; in addition, we are given a vector field $K$ on
$\partial M$ with $L_K\alpha = 0$ and whose flow $(\phi_t^K)$ is one-periodic. One can
always find a collar $j: (-\epsilon;0] \times \partial M \hookrightarrow M$, for some
$\epsilon>0$, such that $j^*\o = d(e^r\alpha)$. Choose a function $H \in \smooth(M,\R)$
with $H(j(r,x)) = e^r(i_K\alpha)(x)$ for all $r \geq -\epsilon/2$. Then the Hamiltonian
flow $(\phi_t^H)$ satisfies $\phi_t^H(j(r,x)) = j(r,\phi_t^K(x))$, and in particular
$\phi_{\pm 1}^H(r,x) = (r,x)$, for all $r \geq -\epsilon/2$. Set
\begin{equation} \label{eq:chi-k}
\chi_K = \phi_{-1}^H \in \Aut(M,\partial M,\o).
\end{equation}
It is easy to see that the class $[\chi_K] \in \pi_0(\Aut(M,\partial M,\o))$ is
independent of the choice of $H$ and $j$. Note that by taking a suitable choice, one
can achieve that $\chi_K$ is the identity outside an arbitrarily small neighbourhood of
$\partial M$. The question we are interested in is: {\em when is $[\chi_K]$ nontrivial,
and more generally, what is its order in $\pi_0(\Aut(M,\partial M,\o))$?} Answering
this can be easy or difficult, depending on the specific situation. We list a few easy
cases:

\begin{examples} \label{ex:chi}
\begin{remarkslist}
\item \label{item:ball}
$[\chi_K]$ is trivial whenever $(\phi_t^K)$ can be extended to a Hamiltonian circle
action on $M$. The simplest example is when $M$ is the unit ball in $\C^n$ and
$(\phi_t^K)$ is any $\C$-linear circle action on $S^{2n-1}$.

\item \label{item:surface}
Take $M$ to be a compact surface with $(\partial M,\alpha) = (S^1,dt)$. Consider the
standard circle action on $\partial M$. Then $\chi_K$ is a positive Dehn twist along a
curve parallel to $\partial M$. $[\chi_K]$ is trivial if $M$ is a disc; in all other
cases, a topological argument shows that it has infinite order.

\item \label{item:variation}
Any continuous map $f: M \longrightarrow M$ which satisfies $f|\partial M = \id$
induces a variation homomorphism $\var(f): H^*(M) \longrightarrow H^*(M,\partial M)$.
For simplicity, consider a smooth map $f$ and cohomology with real coefficients; then
the variation is defined by $\var(f)\theta = \theta - f^*\theta$ for a differential
form $\theta$. The variation of $f = \chi_K^m$ has a particularly simple expression in
the case when the circle action $(\phi_t^K)$ is free:
\begin{equation} \label{eq:pairing}
\int_M \var(\chi_K^m)(a) \cup b = m \int_{\partial M/S^1} p(a) \cup p(b) \quad \text{
for } a,b \in H^*(M;\R).
\end{equation}
Here $p$ is the map $H^*(M;\R) \longrightarrow H^*(\partial M;\R) \longrightarrow
H^{*-1}(\partial M/S^1;\R)$. Note that if $\var(f)$ is nonzero, $f$ cannot be homotoped
to the identity rel $\partial M$. It follows that $[\chi_K]$ has infinite order
whenever there are classes $a,b \in H^*(M;\R)$ with $\int_{\partial M/S^1} p(a) \cup
p(b) \neq 0$.

\parindent0cm \parskip1em
As an application consider $\CP{n}$ with the Fubini-study form $\o_{FS}$. Let $Q
\subset \CP{n}$ be a smooth complex hypersurface of degree $d$, with normal bundle $L$.
A choice of unitary connection $A$ on $L$ with curvature $(i/2\pi d)F_A = -\o_{FS}$
determines a symplectic form $\o_U$ on a neighbourhood $U = \{\xi \in L \suchthat |\xi|
< \epsilon\}$ of the zero-section. If $\epsilon$ is sufficiently small, there is a
symplectic embedding of $(U,\o_U)$ into $\CP{n}$ which forms a tubular neighbourhood of
$Q$. The complement $M = \CP{n} \setminus j(U)$, with $\o = \o_{FS}|M$, has a natural
structure of symplectic manifold with contact type boundary. Moreover, the circle
action $(\phi_t^K)$ on $\partial M$ coming from the obvious circle action on $L$
preserves the contact form. This is a standard construction; see \cite[Lemma
2.6]{mcduff91b} for details. Under Poincar{\'e} duality, the map $p: H^*(M;\R)
\longrightarrow H^{*-1}(\partial M/S^1;\R)$ defined above corresponds to the boundary
operator $\partial: H_*(\CP{n},Q;\R) \longrightarrow H_{*-1}(Q;\R)$. Using
\eqref{eq:pairing} one concludes that $[\chi_K]$ has infinite order whenever there are
middle-dimensional homology classes $\bar{a}, \bar{b} \in H_{n-1}(Q;\R)$ which satisfy
$\bar{a} \cdot \bar{b} \neq 0$ and, if $n$ is odd, also $\leftsc \o^{(n-1)/2}_{FS},
\bar{a} \rightsc = \leftsc \o^{(n-1)/2}_{FS}, \bar{b} \rightsc = 0$. Such classes exist
for all $d \geq 3$, and also for $d = 2$ when $n$ is odd. For $d = 1$ one gets the unit
ball with the standard circle action, so that $[\chi_K]$ is trivial by \ref{item:ball}.
We will show later, using Floer cohomology and graded Lagrangian submanifolds, that
$[\chi_K]$ has infinite order in the remaining case ($d = 2$ and $n$ even).
\end{remarkslist}
\end{examples}

From now on assume that $\partial M$ is connected, $H^1(M) = 0$, and that $2c_1\mo =
0$. Then there is a unique $\infty$-fold Maslov covering $\LL^\infty$ on $M$. There are
two ways to choose an $\LL^\infty$-grading for $\chi_K$. One way is to use Remark
\ref{th:boundary} which says that there is a unique grading $\tchi_K \in
\Autgr(M,\partial M,\o;\LL^\infty)$. Alternatively one can lift $(\phi_t^H)_{-1 \leq t
\leq 0}$ to an isotopy of $\LL^\infty$-graded symplectic automorphisms $(\tphi_t^H)$
with $\tphi_0^H = \id$. Then $\tphi_{-1}^H$ is again a grading of $\chi_K$. However,
this grading does not necessarily lie in $\Autgr(M,\partial M,\o;\LL^\infty)$. In other
words the action of $\tphi_{-1}^H$ on $\LL^\infty_x$, where $x$ is any point in
$\partial M$, may be a nonzero shift. This explains that the two approaches may lead to
different gradings. Let $\sigma_K \in \Z$ be the difference, that is to say
\begin{equation} \label{eq:interior-shift}
\tchi_K = \tphi_{-1}^H \circ [\sigma_K].
\end{equation}

\begin{lemma} \label{th:self-shift}
Let $L$ be a Lagrangian submanifold of $M$ which admits an $\LL^\infty$-grading $\tL$.
If the class $[\chi_K^m] \in \pi_0(\Aut(M,\partial M,\o))$ is trivial for some $m \geq
1$, then $\tL$ is isotopic to $\tL[m\sigma_K]$ as an $\LL^\infty$-graded Lagrangian
submanifold.
\end{lemma}

\proof Since the statement is independent of the choices made in the definition of
$\chi_K$, we can assume that the embedding $j$ satisfies $L \cap \im(j) = \emptyset$.
This means that $\phi_t^H(L) = L$ for all $t$. Since $\tphi_0^H = \id$, it follows that
$\tphi_t^H(\tL) = \tL$ for all $t$. Because of \eqref{eq:interior-shift} one has
$\tchi_K^m(\tL) = \tL[m\sigma_K]$. By assumption there is an isotopy $(\psi_t)$ in
$\Aut(M,\partial M,\o)$ from $\psi_0 = \id$ to $\psi_1 = \chi_K^m$. One can lift this
to an isotopy $(\tpsi_t)$ of $\LL^\infty$-graded symplectic automorphisms with $\tpsi_0
= \id$. This isotopy will remain inside $\Autgr(M,\partial M,\o;\LL^\infty)$, which
implies that $\tpsi_1 = \tchi_K^m$. Hence $\tpsi_t(\tL)$ is an isotopy of
$\LL^\infty$-graded Lagrangian submanifolds from $\tL$ to $\tchi_K^m(\tL) =
\tL[m\sigma_K]$. \qed

\begin{theorem} \label{th:chi}
Let $(M,\o,\alpha)$ be a compact symplectic manifold with contact type boundary, and
$(\phi_t^K)$ a circle action on $\partial M$ which preserves $\alpha$. Assume that
$\partial M$ is connected, that $H^1(M) = 0$, that $2c_1\mo = 0$, and that $[\o] \in
H^2(M;\R)$ is zero. Furthermore, assume that $M$ contains a Lagrangian submanifold $L$
with $H^1(L) = 0$. Let $\chi_K$ be the automorphism of $M$ defined in \eqref{eq:chi-k},
and $\sigma_K$ the number from \eqref{eq:interior-shift}. If $\sigma_K \neq 0$ then
$[\chi_K] \in \pi_0(\Aut(M,\partial M,\o))$ is an element of infinite order.
\end{theorem}

\proof The assumption $[\o] = 0$ implies that the Floer cohomology $HF(L_0,L_1)$ of any
pair of Lagrangian submanifolds $L_0,L_1$ with $H^1(L_0) = H^1(L_1)$ is well-defined
and invariant under Lagrangian isotopy. For $L_0 = L_1$ one has $HF(L_0,L_1) \iso
H^*(L_0;\Z/2)$. Of course, the same properties are true for graded Floer cohomology.
Choose an $\LL^\infty$-grading $\tL$ of $L$. If $[\chi_K]$ was trivial then by Lemma
\ref{th:self-shift} $\tL$ would be graded Lagrangian isotopic to $\tL[\sigma_K]$, and
hence
\[
HF^*(\tL,\tL) \iso HF^*(\tL,\tL[\sigma_K]) = HF^{*+\sigma_K}(\tL,\tL).
\]
Because $HF^*(\tL,\tL)$ is nonzero and concentrated in finitely many degrees, this
implies that $\sigma_K$ must be zero. Conversely, if $\sigma_K \neq 0$ then $[\chi_K]$
must be nontrivial. Since the same argument works for the iterates $\chi_K^m$, it also
follows that $[\chi_K]$ has infinite order. \qed

Let $C^\infty \in H^1(\LL)$ be the global Maslov class of $\LL^\infty$. Take a point $x
\in \partial M$ and a Lagrangian subspace $\Lambda \in \LL_x$, and define $\lambda: S^1
\longrightarrow \LL$ by $\lambda(t) = D\phi_t^H(\Lambda)$. One can easily show that
\begin{equation} \label{eq:boundary-shift}
\sigma_K = -\leftsc C^\infty, [\lambda] \rightsc.
\end{equation}
This formula is useful for determining $\sigma_K$, which is important in applications
of Theorem \ref{th:chi}.

\begin{remarks}
\begin{remarkslist}
\item
The conditions in Theorem \ref{th:chi} were chosen for their simplicity, and are not
the most general ones. For instance, an inspection of the proof shows that the
assumption $H^1(M) = 0$ can be omitted (however, for $H^1(M) \neq 0$ the shift
$\sigma_K$ may depend on the choice of Maslov covering). The condition $[\o] = 0$ is
there to ensure that there is a well-behaved Floer theory, and can also be weakened
considerably. In contrast, the existence of the Lagrangian submanifold $L$, and the
assumption that $2c_1\mo = 0$, are both essential parts of our argument.

\item
There is an alternative approach which dispenses with Lagrangian submanifolds
altogether, and uses instead the Floer cohomology of the automorphism $\chi_K$. This
approach is more difficult to carry out than the one used here, but it is possibly more
general.
\end{remarkslist}
\end{remarks}

\subsection{Manifolds with periodic geodesic flow\label{subsec:cotangent}}

Let $(N^n,g)$ be a Riemannian manifold such that any geodesic of length one is closed.
Take $M = \ucotangent N = \{ \xi \in \cotangent N \suchthat |\xi|_g \leq 1\}$ with the
standard symplectic form $\o \in \Omega^2(M)$ and contact form $\alpha \in
\Omega^1(\partial M)$. The flow of the Reeb vector field $R$ of $\alpha$ is just the
geodesic flow on $N$, and our assumption is that $\phi_1^R = \id$. Hence one can define
an automorphism $\chi_R \in \Aut(M,\partial M,\o)$. Extend $(\phi_t^R)$ radially to a
Hamiltonian circle action $\sigma$ on $M \setminus N$ (the moment map of $\sigma$ is
$\mu(\xi) = |\xi|$). Choose a function $r \in \smooth([0;1],\R)$ with $r(t) = 0$ for $t
\leq 1/3$, $r(t) = -1$ for $t \geq 2/3$. Then a possible choice for $\chi_R$ is
\begin{equation} \label{eq:explicit}
\chi_R(\xi) = \begin{cases} \sigma_{r(|\xi|)}(\xi) & \xi \notin N,\\ \xi & \xi \in N.
\end{cases}
\end{equation}
$M$ admits a Lagrangian distribution $S \subset TM$, formed by the tangent spaces along
the fibres of the projection $M \longrightarrow N$. As indicated in Example
\ref{ex:lagrangian-distribution}, one can use $S$ to define an $\infty$-fold Maslov
covering $\LL^\infty$. In order to satisfy the conditions of Theorem \ref{th:chi}, we
will assume that $H^1(N)$ is zero, so that there is only one $\infty$-fold Maslov
covering. The global Maslov class $C^\infty$ is represented by the fibrewise Maslov
cycle $I = \bigcup_{\xi} \{ \Lambda \in \LL_{\xi} \suchthat \Lambda \cap S_{\xi} \neq
0\}$, and one can write \eqref{eq:boundary-shift} as
\[
\sigma_R = -[I] \cdot [\lambda].
\]
By definition $\lambda(t) = D\sigma_t(\Lambda)$ for some $\Lambda \in \LL_\xi$, $\xi
\in \partial M$. It is convenient to choose $\Lambda$ in the following way: split
$TM_\xi$ into its horizontal and vertical parts, both of which are naturally isomorphic
to $TN_x$, where $x \in N$ is the base point of $\xi$. Then take $\Lambda =
\{(\eta_1,\eta_2) \in TN_x \times TN_x \suchthat \eta_1 = g(\xi,\eta_2)\xi\}$. This has
the consequence that $\lambda(t) \in I$ iff $c(t)$ and $c(0)$ are conjugate points. It
is a familiar result, see \cite[section 4]{duistermaat76} or \cite[section
6]{robbin-salamon95}, that the local intersection number at a point $\lambda(t) \in I$
is equal to the multiplicity $m_c(t) > 0$ of $c(t)$ as a conjugate point of $c(0)$.
Therefore
\begin{equation} \label{eq:morse}
\sigma_R = -\sum_t m_c(t),
\end{equation}
where the sum is over all conjugate points $t \in S^1$. Since $m_c(0) = n-1$,
$\sigma_R$ is always negative ($N = S^1$ is impossible because we have assumed that
$H^1(N) = 0$). $M$ always contains a Lagrangian submanifold (the zero-section) with
zero first Betti number. Moreover, $[\o] = 0$. By applying Theorem \ref{th:chi} one
obtains

\begin{cor} \label{th:besse}
Let $(N,g)$ be a closed Riemannian manifold with $H^1(N) = 0$ such that any geodesic of
length one is closed. Then, for $M = \ucotangent N$ and the Reeb vector field $R$,
$[\chi_R] \in \pi_0(\Aut(M,\partial M,\o))$ has infinite order. \qed \end{cor}

The main examples of manifolds with periodic geodesic flow are compact globally
symmetric spaces of rank one \cite{besse}. These spaces are two-point homogeneous,
which means that the isometry group $\mathrm{Iso}(N,g)$ acts transitively on $\partial
M$. Hence all geodesics have the same minimal period. Symmetric spaces also have the
property that any geodesic path $c: [0;T] \longrightarrow N$ with $c(T) = c(0)$ is a
closed geodesic \cite[p. 144]{klingenberg}. It follows that the energy functional on
the {\em based} loop space $\Omega N$ is a Morse-Bott function. Its critical point set
consists of one point (the constant path) and infinitely many copies of $S^{n-1}$. The
Morse index of the point is zero. The Morse indices of the $(n-1)$-spheres can be
computed by comparing \eqref{eq:morse} with Morse's index theorem: they are
$-\sigma_R-n+1$, $-2\sigma_R-n+1$, etc. This means that for $n \geq 3$, $\Omega N$ has
a CW-decomposition with one cell of dimension $-\sigma_R-n+1$ and other cells of
dimension $\geq -\sigma_R-n+3$. It follows that $\Omega N$ is precisely
$(-\sigma_R-n)$-connected, and that $N$ is precisely $(-\sigma_R-n+1)$-connected. This
approach, complemented by explicit computations for $\RP{2}$ and $S^2$, yields the
following values for $\sigma_R$:
\[
\begin{tabular}{r||c|c|c|c|c}
$N$ & $S^m$ & $\RP{m}$ & $\CP{m}$ & $\HP{m}$ & $F_4/\mathrm{Spin}_9$ \\ \hline
$\sigma_R$ & $2-2m$ & $1-m$ & $-2m$ & $-4m+2$ & $-22$
\end{tabular}
\]
Here $F_4/\mathrm{Spin}_9$ is the exceptional symmetric space diffeomorphic to the
Cayley plane (which is 16-dimensional and 7-connected). For $N = S^{2m+1}$ or
$\RP{2m+1}$ the result of Corollary \ref{th:besse} can be obtained more easily, without
using any symplectic geometry, by computing the variation $\var(\chi_R)$. For $N =
\RP{2m}$ the space $M$ can be identified with the complement of a neighbourhood of a
quadric hypersurface $Q \subset \CP{2m}$, thus settling the remaining case in Example
\ref{ex:chi}\ref{item:variation}. The most interesting example is that of $\CP{m}$,
since there we can show that the non-vanishing of $[\chi_R]$ is a genuinely symplectic
phenomenon:

\begin{prop} \label{th:topologically-trivial}
Let $M = \ucotangent \CP{m}$, $m \geq 1$. Then $\chi_R$ can be deformed to the identity
in the group $\Diff(M,\partial M)$ of diffeomorphisms which act trivially on $\partial
M$.
\end{prop}

\begin{lemma} \label{th:family-of-circle-actions} Let $M = \ucotangent \CP{m}$. Then
there is a smooth family $\sigma^s$, $0 \leq s \leq 1$, of circle actions on $\partial
M$ with $\sigma^0_t = \phi_t^R$, and such that $\sigma^1$ extends to a circle action on
the whole of $M$.
\end{lemma}

\proof Since $\CP{m} = S^{2m+1}/S^1$, $\cotangent\CP{m}$ is a symplectic quotient
$\cotangent S^{2m+1}/\!/S^1$. Explicitly $M = \{ (u,v) \in \C^{m+1} \times \C^{m+1}
\suchthat |u| = 1, \; |v| \leq 1, \; \leftsc u,v \rightsc_{\C} = 0 \} / S^1$, where
$S^1$ acts by $t \cdot (u,v) = (e^{2\pi i t}u, e^{2\pi i t}v)$. The diagonal
$SU(2)$-action on $\C^{m+1} \times \C^{m+1}$ descends to an action of $PU(2) =
SU(2)/\!\pm\!1$ on $\partial M$, which we will denote by $\rho$. Assume that we have
rescaled the standard symplectic form to make the Reeb flow one-periodic, that is, $\o
= \pi/2 \sum_j (dv_j \wedge d\bar{u}_j + d\bar{v}_j \wedge du_j)$. Then the Reeb flow
is $\phi_t^R = \rho(\exp(t A))$ with $A = \left(\begin{smallmatrix} 0 & -\pi \\ \pi & 0
\end{smallmatrix}\right) \in su_2$. The family $\sigma^s$ of circle actions is defined
by $\sigma^s_t = \rho(\exp(t A^s))$, where $(A^s)$ is any path in $su(2)$ from $A^0 =
A$ to $A^1 = \mathrm{diag}(\pi i,-\pi i)$ such that $\exp(A^s) = -1$ for all $s$.
Because $\partial M$ is a quotient one can write $\sigma^1_t(u,v) = (u,e^{-2\pi i
t}v)$, and this shows that $\sigma^1$ extends to a circle action on all of $M$. \qed

\proof[Proof of Proposition \ref{th:topologically-trivial}] Clearly, using the family
$\sigma^s$ of circle actions one can deform $\chi_R$ inside $\Diff(M,\partial M)$ to
$\chi_R'(\xi) = \sigma^1_{r(|\xi|)}(\xi)$. To deform this to the identity one uses the
isotopy $\psi_s(\xi) = \sigma^1_{(1-s)r(|\xi|)-s}(\xi)$. \qed

\begin{remark}
It seems likely that the maps $\chi_R$ on $\ucotangent\CP{m}$ are not only
differentiably isotopic to the identity but also `fragile' in the sense of
\cite{seidel97}. We have not checked the details.
\end{remark}

\subsection{Weighted homogeneous polynomials\label{subsec:weighted}}

A polynomial $p \in \C[x_0,\dots,x_n]$, $n \geq 1$, is called weighted homogeneous if
there are integers $\beta_0,\dots,\beta_n,\beta > 0$ such that $p(z^{\beta_0} x_0,
\dots, z^{\beta_n} x_n) = z^\beta p(x_0,\dots,x_n)$. The numbers $w_i = \beta_i/\beta$
are called the weights of $p$. Throughout this section $p(x)$ will be a weighted
homogeneous polynomial with an isolated critical point at $x = 0$; because of the
homogeneity, this implies that $p$ has no other critical points. By definition, the
link of the singular point is $L = p^{-1}(0) \cap S^{2n+1}$. This is in fact a contact
submanifold of $S^{2n+1}$ with respect to the standard contact form $\alpha_{S^{2n+1}}
= \frac{i}{4} \sum_j z_jd\bar{z}_j - \bar{z}_j dz_j$. Let $B^{2n+2} \subset \C^{n+1}$
be the closed unit ball. Fix a cutoff function $\psi$ with $\psi(t^2) = 1$ for $t \leq
1/3$ and $\psi(t^2) = 0$ for $t \geq 2/3$. For $z \in \C \setminus \{0\}$ set $F_z = \{
x \in B^{2n+2} \suchthat p(x) = \psi(|x|^2)z \}$.

\begin{lemma} \label{th:milnor-fibre}
There is an $\epsilon>0$ such that for all $0 < |z| < \epsilon$, $F_z$ is a symplectic
submanifold of $B^{2n+2}$ with boundary $\partial F_z = L$.
\end{lemma}

\proof The Zariski tangent space of $F_z$ is $(TF_z)_x = \ker L(x,z)$ with $L(x,z):
\C^{n+1} \longrightarrow \C$, $L(x,z)\xi = dp(x)\xi - 2z\,\psi'(|x|^2)\leftsc x,\xi
\rightsc$. Clearly $L(x,z)$ is onto in the following three cases: (1) $z = 0$ and $x
\neq 0$, (2) $0 < |x| \leq 1/3$, (3) $|x| \geq 2/3$. On the other hand, the set of
those $(x,z)$ for which $L(x,z)$ is onto must be open. This implies that $L(x,z)$ is
onto for all $(x,z) \in F_z$, provided that $z \neq 0$ is sufficiently small. An
argument of the same kind shows that $(TF_z)_x$ is a symplectic subspace of $\C^{n+1}$
for all small $z \neq 0$. \qed

We fix a $z_0$ with $0 < |z_0| < \epsilon$ and call $M = F_{z_0}$ the Milnor fibre of
the singular point $0 \in p^{-1}(0)$. Clearly $(M,\o = \o_{\C^{2n+2}}|M,\alpha =
\alpha_{S^{2n+1}}|L)$ is a compact symplectic manifold with contact type boundary. The
choice of $z_0$ is not really important since one can prove (by a standard argument
using Moser's technique) that any two choices give symplectically isomorphic Milnor
fibres. During the following discussion, we will repeatedly make use of our right to
pass to a smaller $z_0$. The next Lemma shows that $M$ is diffeomorphic to what is
traditionally called the Milnor fibre.

\begin{lemma} \label{th:classical-milnor-fibre}
There is an $\epsilon>0$ such that for all $0 < |z| < \epsilon$, $F_z$ is diffeomorphic
to $p^{-1}(z) \cap B^{2n+2}$.
\end{lemma}

\proof For $(z,t) \in (\C \setminus \{0\}) \times [0;1]$ consider $G_{(z,t)} = \{ x \in
B^{2n+2} \suchthat p(x) = t\psi(|x|^2)z + (1-t)z\}$. Using the same argument as before,
one can prove that these are smooth manifolds for all sufficiently small $z$. If we fix
such a $z$, the $G_{(z,t)}$ form a differentiable fibre bundle over $[0;1]$. Hence
$G_{(z,1)} = F_z$ and $G_{(z,0)} = p^{-1}(z) \cap B^{2n+2}$ are diffeomorphic. \qed

\begin{cor} \label{th:milnor}
$M$ is $(n-1)$-connected. In fact, it is homotopy equivalent to a nontrivial wedge of
$n$-spheres.
\end{cor}

This follows from Lemma \ref{th:classical-milnor-fibre} and a classical result of
Milnor \cite[Theorem 6.5]{milnor68}. Denote by $A = \{x \in M \suchthat 2/3 \leq |x|
\leq 1\}$ the part of $M$ on which it agrees with $p^{-1}(0)$. Let $\theta_{\C^{n+1}}$
be the complex $n$-form on $\C^{n+1}$ given by $(\theta_{\C^{n+1}})_x(v_1,\dots,v_n) =
\det_{\C}(\overline{dp}_x \wedge v_1 \wedge \dots \wedge v_n)$.

\begin{lemma} \label{th:theta-form}
There is an $\o$-compatible almost complex structure $J$ on $M$, and a nowhere
vanishing $J$-complex $n$-form $\theta$, which have the following properties: the
restriction of $J$ to $A$ agrees with the standard complex structure on $A \subset
p^{-1}(0)$, and $\theta|A = \theta_{\C^{n+1}}|A$.
\end{lemma}

\proof Consider the $2n$-dimensional complex vector bundle $K \longrightarrow M$ with
fibres $K_x = \ker(dp_x)$. $\theta_{\C^{n+1}}$ defines a complex $n$-form on every
fibre of $K$, and these $n$-forms are easily seen to be nonzero. Both $K$ and $TM$ are
subbundles of the trivial bundle $\C^{n+1} \times M$. The proof of Lemma
\ref{th:milnor-fibre} shows that, by choosing $z_0$ sufficiently small, one can make
these two subbundles arbitrarily close. This means that the orthogonal projection $P_K:
TM \longrightarrow K$ is a bundle isomorphism, and that the pullback of the complex
structure on $K$ is an almost complex structure $J'$ on $TM$ which is $\o$-tame. The
pullback $\theta' = P_K^*(\theta_{\C^{n+1}}|K)$ is a $J'$-complex $n$-form on $M$ which
is nowhere zero. Since $P_K$ is the identity over $A$, $J'$ and $\theta'$ have all the
properties required in the Lemma, except that $J'$ may not be $\o$-compatible. However,
one can easily find an $\o$-compatible almost complex structure $J$ which agrees with
$J'$ on $A$. Moreover, $J$ and $J'$ can be deformed into each other through almost
complex structures, and the deformation can be chosen constant on $A$. This implies
that there is a bundle automorphism $Q: TM \longrightarrow TM$ such that $Q^*(J') = J$,
and which is the identity on $A$. Now set $\theta = Q^*(\theta')$. \qed

\begin{cor} \label{th:chern}
$c_1(M,\o) = 0$.
\end{cor}

This is clear, since there is a nowhere vanishing complex $n$-form.

\begin{lemma} \label{th:vanishing-cycle}
$(M,\o)$ always contains an embedded Lagrangian $n$-sphere.
\end{lemma}

\proof[Sketch of proof.] The method which produces such spheres is to deform $p$ by
adding a linear term $\lambda$, such that $p + \lambda$ has only nondegenerate critical
points. The Lagrangian spheres appear as vanishing cycles associated to these critical
points. We will now explain one version of this argument; for variations on this theme
see \cite[section 1.4]{seidel97}. Choose some small $\delta > 0$, and let $D = \{
(z,\lambda) \in \C \times (\C^{n+1})^* \suchthat |z| < \delta, |\lambda| < \delta\}$.
For $(z,\lambda) \in D$ define
\[
\tilde{F}_{(z,\lambda)} = \{ x \in B^{2n+2} \suchthat p(x) + \lambda(x) = \psi(|x|^2)
z\}.
\]
If $\delta$ is sufficiently small, there are two possibilities for each $(z,\lambda)$:
either (1) $\tilde{F}_{(z,\lambda)}$ is a smooth symplectic submanifold of $B^{2n+2}$,
or (2) the complex hypersurface $p(x) + \lambda(x) = z$ has a singular point $x$ with
$|x| \leq 1/3$. The subset $\Delta \subset D$ where (2) occurs is a complex
hypersurface. Hence $D \setminus \Delta$ is connected. By an application of Moser's
technique, it follows that all the symplectic manifolds $\tilde{F}_{(z,\lambda)}$
occurring in case (1) are isomorphic to the Milnor fibre $M = \tilde{F}_{(z_0,0)}$.

Take a generic small $\lambda$. Then $p + \lambda$ has only nondegenerate (Morse-type)
critical points. These critical points lie close to the origin, and there is always at
least one of them (this is a well-known fact, which follows from considering the Milnor
number $\mu$ of the singularity). Choose a critical point $x$ of $p+\lambda$, such that
$|x| < 1/4$, and set $z = p(x) + \lambda(x)$. Since $x$ is non-degenerate, one can
write
\[
(p+\lambda)(x+y) = z + Q(y) + \text{(higher order terms in $y$)}
\]
where $Q$ is a nondegenerate complex quadratic form. A careful application of Moser's
technique shows that for all $0 < \epsilon_2 \ll \epsilon_1 \ll |\lambda|$, one can
embed $U = \{ y \in \C^{n+1} \suchthat |y| \leq \epsilon_1, \; Q(y) = \epsilon_2\}$
symplectically into $\tilde{F}_{(z+\epsilon_2,\lambda)}$. Now $U$ is symplectically
isomorphic to a neighbourhood of the zero-section in $\cotangent\!S^n$, hence contains
a Lagrangian $n$-sphere. It follows that $\tilde{F}_{(z+\epsilon_2,\lambda)}$, and
hence the Milnor fibre, contain a Lagrangian $n$-sphere. \qed

The Milnor fibration associated to the singular point $0 \in p^{-1}(0)$ is obtained by
putting together the manifolds $F_z$ for all $|z| = |z_0|$:
\[
p: F = \bigcup_{|z| = |z_0|} F_z \times \{z\} \longrightarrow |z_0| S^1.
\]
The proof of Lemma \ref{th:milnor-fibre} shows that this is a smooth and proper
fibration. Moreover, if we pull back $\o_{\C^{n+1}}$ to $F$ via the obvious projection,
we obtain a closed two-form $\Omega$ whose restriction to any fibre is a symplectic
form. Such a two-form defines a connection $TF^h \subset TF$, given by the
$\O$-orthogonal complements of the tangent spaces along the fibres:
\[
TF^h_{(x,z)} = \{X \in TF_{(x,z)} \suchthat \O(X,Y) = 0 \text{ for all $Y$ such that }
Dp(Y) = 0 \}.
\]
The fact that $\O$ is closed implies that the parallel transport maps of this
connection are symplectic isomorphisms between the fibres. In our case, since the
fibres are manifolds with boundary, one has to check that the parallel transport maps
exist. But this it clear because near the boundary $\partial F = \bigcup_z \partial F_z
\times \{z\}$ one has a natural trivialization $A \times |z_0|\,S^1 \subset F$, and the
connection is compatible with this trivialization. This also implies that if we go once
around the base $|z_0| S^1$, the parallel transport yields a map $f \in \Aut(M,\partial
M,\o)$. We call $f$ the {\em symplectic monodromy} of our singularity. One can show
that the class $[f] \in \pi_0(\Aut(M,\partial M,\o))$ is, in a suitable sense,
independent of the choices of $\psi$ and $z_0$.

Up to now, the fact that $p$ is weighted homogeneous has not been of any importance. We
will now begin to exploit this particular feature of our situation. Let $\sigma$ be the
complex circle action on $\C^{n+1}$ with multiplicities $\beta_0,\dots, \beta_n$, and
$K$ the vector field generating it. $\sigma|L$ is a circle action preserving $\alpha$.
Hence one can construct the associated automorphism $\chi_K \in \Aut(M,\partial M,\o)$
of the Milnor fibre. We will use a particular choice which is $\chi_K = \phi_{-1}^H$,
where $H \in \smooth(M,\R)$ is given by $H(x) = \pi \sum_i \beta_i |x_i|^2$. From now
on assume that $n \geq 2$. Since $H^1(M)$ is zero (Corollary \ref{th:milnor}),
$\partial M$ is connected (this follows from Corollary \ref{th:milnor} by a
Poincar{\'e} duality consideration), and $2c_1\mo = 0$ (Corollary \ref{th:chern}), one
can define the shift $\sigma_K$ of $\chi_K$.

\begin{lemma} \label{th:weights}
$\sigma_K = 2(\beta - \sum_j \beta_j)$.
\end{lemma}

\proof Let $J,\theta$ be as in Lemma \ref{th:theta-form}. As explained in Example
\ref{ex:calabi-yau}, $\Theta = \theta^{\otimes 2}/|\theta|^2$ defines an $\infty$-fold
Maslov covering on $M$, whose global Maslov class is represented by the map
$\det_\Theta^2: \LL \longrightarrow S^1$. Together with \eqref{eq:boundary-shift} this
means that $-\sigma_K$ is the degree of the map $c: S^1 \longrightarrow S^1$, $c(t) =
\det_\Theta^2(D\phi_t^H(\Lambda))$ for some $\Lambda \in \LL_x$, $x \in \partial M
\subset A$. The restriction of $\phi_t^H$ to $A$ is simply the circle action $\sigma$,
and $\theta|A$ agrees with $\theta_{\C^{n+1}}$. An explicit computation shows that
$\sigma_t^*(\theta_{\C^{n+1}}) = e^{2\pi i t(\beta_0 + \dots + \beta_n -
\beta)}\theta_{\C^{n+1}}$. Therefore $c$ has the form $c(t) = e^{4\pi i t(\beta_0 +
\dots + \beta_n - \beta)}c(0)$. \qed

\begin{lemma} \label{th:the-monodromy} $\chi_K$ is equal to the $\beta$-th iterate $f^\beta$ of
the symplectic monodromy.
\end{lemma}

\proof Let $\tilde{X} \in \smooth(TE^h)$ be the unique horizontal lift of the vector
field $X(z) = 2\pi i z$ on $|z_0|S^1$. Its flow $(\mu_t)$ maps $F_z$ to $F_{e^{2\pi i
t}z}$ symplectically. By definition, the symplectic monodromy is $f = \mu_1|F_{z_0}$.
Now let $\rho$ be the circle action on $F$ given by $\rho_t(x,z) = (\sigma_t(x),
e^{2\pi i \beta t}z)$. We denote the Killing vector field of $\rho$ by $Y$. Since
$\rho$ preserves $\O$, the connection $TE^h$ is $\rho$-invariant. It follows that
$\mu_t$ commutes with $\rho$ for any $t$. Therefore $\eta_t = \rho_{-t} \circ
\mu_{\beta t}$ is the flow on $F$ generated by $\beta\tilde{X}-Y$. Note that $\eta_t$
maps any fibre $F_z$ to itself symplectically. Let $H$ be the function which we have
used to define $\chi_K$. Clearly $d(H|F_{z_0}) = -i_Y\O|F_{z_0}$. Since
$i_{\tilde{X}}\O$ vanishes on each fibre $F_z$ one also has
\[
d(-H|F_{z_0}) = (-i_{\beta\tilde{X}-Y}\O)|F_{z_0}.
\]
This means that $(\eta_t|F_{z_0})$ is the Hamiltonian flow of $-H|F_{z_0}$. Hence by
definition $\chi_K = \eta_1|F_{z_0}$. On the other hand, by the definition of $\eta_t$,
$f^\beta = \eta_1|F_{z_0}$. \qed

We can now prove our main result about symplectic monodromy.

\begin{theorem} \label{th:monodromy}
Let $p \in \C[x_0,\dots,x_n]$ be a weighted homogeneous polynomial with an isolated
critical point at $0$. Assume that $n \geq 2$ and that the sum of the weights $w_i$ is
not one. Then the symplectic monodromy $f$ defines a class of infinite order in
$\pi_0(\Aut(M,\partial M,\o))$.
\end{theorem}

\proof Since $f^\beta = \chi_K$, it is sufficient to prove that $[\chi_K]$ has infinite
order. Lemma \ref{th:weights} shows that $\sigma_K = 2\beta(1-\sum_i w_i) \neq 0$.
Therefore one only needs to apply Theorem \ref{th:chi}. The necessary assumptions about
$M$ have all been proved above except for $[\o] = 0$, which is obvious from the
definition. \qed

\begin{examples-remarks}
\begin{remarkslist} \item
Let $p \in \C[x_0,x_1,x_2]$ be one of the standard models for the du Val (or simple)
singularities \cite{arnold72}. These models are weighted homogeneous, and the sum of
the weights is $>1$. For example, if $p(x) = x_0^2 + x_1^3 + x_1x_2^3$ is the
singularity of type $(E_7)$ then $w_0 + w_1 + w_2 = 1/2 + 1/3 + 2/9 = 19/18$. Hence
Theorem \ref{th:monodromy} applies, showing that $[f] \in \pi_0(\Aut(M,\partial M,\o))$
has infinite order. In contrast, it follows from Brieskorn's simultaneous resolution
\cite{brieskorn68} that the class of $f$ in $\pi_0(\Diff(M,\partial M))$ has {\em
finite order}. Hence, at least in these cases, Theorem \ref{th:monodromy} expresses a
genuinely symplectic phenomenon.

\item
It is possible that the assumption $\sum_i w_i \neq 1$ might be removed. In fact, there
are many cases when $\sum_i w_i = 1$ and in which the monodromy has infinite order for
topological reasons, for example $p(x_0,x_1,x_2,x_3) = x_0^4 + x_1^4 + x_2^4 + x_3^4$.

\item
Let $p \in \C[x_0,x_1]$ be a weighted homogeneous polynomial with an isolated critical
point. Then $M$ is a connected surface and not a disc. The iterate $f^\beta$ of the
monodromy can be written as a composition of positive Dehn twists along the connected
components of $\partial M$. Using this one can show easily that $[f] \in
\pi_0(\Diff(M,\partial M))$ is always of infinite order. This means that Theorem
\ref{th:monodromy} holds also for $n = 1$.
\end{remarkslist}
\end{examples-remarks}

\section{Knotted Lagrangian spheres}

\subsection{Generalized Dehn twists\label{sec:dehn-twists}}

This section contains the basic definitions and some facts, both topological and
symplectic, which are used later on. Throughout $\mo$ will be a compact symplectic
manifold of dimension $2n$. By a Lagrangian sphere in $M$ we will mean a Lagrangian
embedding $S^n \hookrightarrow M$. Such embeddings will be denoted by the letters
$l,l_1,l_2,\dots$ and their images by $L,L_1,L_2,\dots$. An {\em
$(A_k)$-configuration}, $k \geq 2$, is a collection of Lagrangian spheres
$(l_1,\dots,l_k)$ which are pairwise transverse, such that $L_i \cap L_j = \emptyset$
for $|i-j| \geq 2$ and $|L_i \cap L_{i{\pm}1}| = 1$. The name comes from the
relationship with singularity theory. In fact the Milnor fibre of the
$(A_k)$-singularity, which is the hypersurface
\begin{equation} \label{eq:ak}
M = \{x \in \C^{n+1} \suchthat x_0^{k+1} + x_1^2 + \dots + x_n^2 = \epsilon, |x| \leq
1\}
\end{equation}
for sufficiently small $\epsilon \neq 0$, contains such a configuration. This was
proved in \cite[Proposition 8.1]{seidel98b} for $n = 2$, and the general case can be
treated in the same way (here we have used the classical form of the Milnor fibre,
rather than the definition adopted in section \ref{subsec:weighted}; this does not
really matter, since the Milnor fibre as defined there also contains an
$(A_k)$-configuration).

Consider $U = \ucotangent S^n = \{ \xi \in \cotangent S^n \suchthat |\xi| \leq 1 \}$
with its standard symplectic structure $\o_U$. The complement of the zero-section $S^n
\subset U$ carries a Hamiltonian circle action $\sigma$ whose moment map is the length
$|\xi|$. In the coordinates
$
U = \{ (u,v) \in \R^{n+1} \times \R^{n+1} \suchthat |v| = 1, \; |u| \leq 1, \;
\gen{u,v} = 0\}
$
with $\o_U = \sum_j du_j \wedge dv_j$ one has $\sigma_t(u,v) = (\cos(2\pi t)u -
\sin(2\pi t)v|u|, \cos(2\pi t)v + \sin(2\pi t)u/|u|)$. Note that $\sigma_{1/2}(u,v) =
(-u,-v)$ extends smoothly over the zero-section. Therefore one can define a
diffeomorphism $\tau$ of $U$ by setting
\[
\tau(\xi) = \sigma_{\psi(|\xi|)}(\xi),
\]
where $\psi$ is a function with $\psi(t) = 1/2$ for $t \leq 1/3$ and $\psi(t) = 0$ for
$t \geq 2/3$. $\tau$ is equal to the identity near $\partial U$, and it acts on the
zero-section as the antipodal map. An explicit computation shows that $\tau$ is
symplectic. We call it a {\em model generalized Dehn twist.} Now let $l$ be any
Lagrangian sphere in $\mo$. One can always find an embedding $j: U \longrightarrow M$
such that $j|S^n = l$ and $j^*\o = \delta\,\o_U$ for some $\delta>0$. By extending $j
\tau j^{-1}$ trivially over $M \setminus \im(j)$ one defines a symplectic automorphism
$\tau_l$ of $M$, which we call a {\em generalized Dehn twist along $l$}. It is not
difficult to see that the class $[\tau_l]$ in $\pi_0(\Aut(M,\o))$ (or in
$\pi_0(\Aut(M,\partial M,\o))$, if $M$ has a boundary) is independent of the choice of
$j$ and $\psi$. For this reason, we will often speak of $\tau_l$ as {\em the}
generalized Dehn twist along $l$. For $n = 1$ these maps are just the ordinary
(positive) Dehn twists along curves on a surface.

\begin{remark} It is an open question whether $[\tau_l]$ depends only on the image $L$. If
$l$ and $l'$ are two embeddings with the same image, and such that $l^{-1} \circ l'$ is
isotopic to the identity in $\Diff(S^n)$, one can easily prove that $[\tau_l] =
[\tau_{l'}]$. Moreover, the same holds if $l^{-1} \circ l'$ is an element of $O(n+1)$
(this shows that, just as in the case $n = 1$, the choice of orientation of $L$ is not
important). For $n \leq 3$ it is known that $\pi_0(\Diff^+(S^n)) = 1$ \cite{munkres58}
\cite{cerf68} which implies that $[\tau_l]$ does indeed depend only on $L$, but in
higher dimensions $\pi_0(\Diff^+(S^n))$ is often nonzero.
\end{remark}

We will first look at generalized Dehn twists from a topological point of view. Since
these maps are a symplectic form of the classical Picard-Lefschetz transformations,
their action on homology is given by the familiar formula
\begin{equation} \label{eq:picard-lefschetz}
(\tau_l)_*(c) = c - (-1)^{n(n-1)/2} (c \cdot [L]) [L].
\end{equation}
Using the fact that $[L] \cdot [L] = (-1)^{n(n-1)/2} \chi(L)$ for any $n$-dimensional
Lagrangian submanifold, one obtains the following

\begin{lemma} \label{th:homology}
\begin{theoremlist} \item \label{item:homology-even}
If $n$ is even then $(\tau_l)_*$ has order two. If $n$ is odd then $(\tau_l)_*$ has
infinite order iff $[L]$ is not a torsion class (otherwise $(\tau_l)_* = \id$).

\item \label{item:homology-odd}
Assume that $n$ is odd, and that $l_1,l_2$ are two Lagrangian spheres with $[L_1] \cdot
[L_2] = \pm 1$. Set $g = (\tau_{l_2}\tau_{l_1})^3$. Then $g^2$ induces the identity on
homology. \qed
\end{theoremlist}
\end{lemma}

For $n = 2$ it is known \cite[Lemma 6.3]{seidel98b} that $\tau_l^2$ is actually
isotopic to the identity in $\Diff(M)$. It seems natural to ask whether this holds for
other even $n$; there is also the analogous problem for the map $g^2$ defined in
\ref{item:homology-odd}. In both cases the answer is unknown to the author. However,
there are some weaker topological results which are easier to obtain, and which are
sufficient for our purpose.

\begin{lemma} \label{th:haefliger}
\begin{theoremlist}
\item \label{item:even}
Assume that $n$ is even, and that there is an $(A_2)$-configuration of Lagrangian
spheres $(l_1,l_2)$ in $M$. Then $l_1^{(k)} = \tau_{l_2}^{2k}\circ l_1$ is isotopic to
$l_1$ through smooth embeddings $S^n \hookrightarrow M$ for any $k \in \Z$.

\item \label{item:odd}
Assume that $n$ is odd and $\geq 5$, and that there is an $(A_3)$-configuration of
Lagrangian spheres $(l_1,l_2,l_3)$ in $M$. Then $l_1^{(k)} = g^{2k}\circ l_1$, where $g
= (\tau_{l_2}\tau_{l_3})^3$, is isotopic to $l_1$ through smooth embeddings $S^n
\hookrightarrow M$ for any $k \in \Z$.
\end{theoremlist}
\end{lemma}

\proof \ref{item:even} If $n = 2$ then $\tau_{l_2}^2$ is isotopic to the identity in
$\Diff(M)$, which implies our result. Hence we can assume that $n \geq 4$. let $W
\subset M$ be a regular neighbourhood of $L_1 \cup L_2$. Since $W$ retracts onto $L_1
\cup L_2$, it is $(n-1)$-connected with $H_n(W) = \Z\gen{L_1} \oplus \Z\gen{L_2}$. We
can assume that $\tau_{l_2}$ has been chosen in such a way that it preserves $W$. Lemma
\ref{th:homology}\ref{item:homology-even} shows that all the embeddings $l_1^{(k)}$
represent the same homology class in $W$. Hence, by Hurewicz's theorem, they are
homotopic as continuous maps $S^n \longrightarrow W$. The proof is completed by
applying a result of Haefliger \cite{haefliger62} which shows that any two homotopic
embeddings $S^n \longrightarrow W$ are differentiably isotopic.

\ref{item:odd} is proved in the same way. The condition $n \geq 5$ is necessary in
order to use Haefliger's result. The result is easily seen to be false for $n = 1$, but
the author was unable to decide the remaining case $n = 3$. In that dimension, there
are obstructions for two homotopic embeddings $S^3 \longrightarrow W$ to be
differentiably isotopic. These obstructions are completely understood in principle, see
\cite[Corollary B]{haefliger66}, but not so easy to compute in practice. \qed

Now consider an $(A_2)$-configuration $(l_1,l_2)$ of Lagrangian spheres in $M$. From
\eqref{eq:picard-lefschetz} it follows that, in any dimension, $[\tau_{l_2}(L_1)] = \pm
[\tau_{l_1}^{-1}(L_2)] \in H_n(M)$. In fact the following stronger result is true:

\begin{lemma} \label{th:pre-braiding}
$\tau_{l_1}(L_2)$ and $\tau_{l_2}^{-1}(L_1)$ are isotopic as (unoriented) Lagrangian
submanifolds of $\mo$. In fact, both of them are Lagrangian isotopic to the surgery
$L_1 \# L_2$. \end{lemma}

This was proved in \cite[Appendix]{seidel98b} for $n = 2$, and the argument given there
carries over to arbitrary $n$. For future use we need to recall one aspect of the
proof: both $\tau_{l_1}(L_2)$ and $L_1 \# L_2$ agree with $L_2$ away from a
neighbourhood of $L_1$ which, by an appropriate choice of $\tau_{l_1}$ and of the
surgery, can be made arbitrarily small. The Lagrangian isotopy constructed in
\cite{seidel98b} between them remains constant outside this neighbourhood. Similarly,
the isotopy between $L_1 \# L_2$ and $\tau_{l_2}^{-1}(L_1)$ is constant outside a
neighbourhood of $L_2$.

\begin{lemma} \label{th:second-pre-braiding}
Assume that $(l_1,l_2)$ form an $(A_2)$-configuration of Lagrangian spheres in $M$, and
set $g = (\tau_{l_1}\tau_{l_2})^3$. Then $g(L_1)$ is Lagrangian isotopic to $L_1$, and
$g(L_2)$ is Lagrangian isotopic to $L_2$.
\end{lemma}

\proof Using Lemma \ref{th:pre-braiding} and the obvious fact that $\tau_l(L) = L$ for
any $l$, one sees that $(\tau_{l_1}\tau_{l_2})^3(L_2) \isotopic
(\tau_{l_1}\tau_{l_2})^2\tau_{l_2}^{-1}(L_1) \isotopic
\tau_{l_1}\tau_{l_2}\tau_{l_1}(L_1) \isotopic \tau_{l_1}\tau_{l_1}^{-1}(L_2) = L_2$,
where $\isotopic$ stands for Lagrangian isotopy. The proof for $L_1$ is similar. \qed
\subsection{The graded point of view\label{subsec:graded-dehn}}

From now on we assume that $(M,\o)$ satisfies $2c_1(M,\o) = 0$ and that its dimension
$2n$ is at least four (the case of classical Dehn twists, $n = 1$, is more complicated
because $\cotangent\!S^1$ admits infinitely many different Maslov coverings). Choose an
$\infty$-fold Maslov covering $\LL^\infty$ on $M$. Lemma \ref{th:newbasic-two} implies
that any Lagrangian submanifold $L \subset M$ with $H^1(L) = 0$ admits an
$\LL^\infty$-grading. Let $l$ be a Lagrangian sphere in $M$ and $\tau_l$ the
generalized Dehn twist along it defined using some embedding $j: U \hookrightarrow M$.
By definition $\tau_l$ is the identity outside $\im(j)$.

\begin{lemma} \label{th:preferred-grading}
There is a unique $\LL^\infty$-grading $\ttau_l$ of $\tau_l$ which acts trivially on
the part of $\LL^\infty$ which lies over $M \setminus \im(j)$.
\end{lemma}

\proof The uniqueness is obvious. To prove the existence, consider the local model $U =
\ucotangent\!S^n$. Since $c_1(U,\o_U) = 0$ and $H^1(U) = 0$, $U$ admits a unique
$\infty$-fold Maslov covering $\LL^\infty_U$. Remark \ref{th:boundary} says that the
model generalized Dehn twist $\tau$ has a unique $\LL_U^\infty$-grading $\ttau$ which
acts trivially on the part of $\LL_U^\infty$ which lies over $\partial U$. Now, for an
arbitrary Lagrangian sphere $l$ in $M$ and Maslov covering $\LL^\infty$, one can
identify $\LL^\infty|\im(j)$ with $\LL^\infty_U$ and then extend $\ttau$ by the
identity to an $\LL^\infty$-grading $\ttau_l$ of $\tau_l$. \qed

We emphasize that this preferred grading $\ttau_l$ does not depend on the choice of a
grading for $L$.

\begin{lemma} \label{th:dehn-shift}
The preferred $\LL^\infty$-grading $\ttau_l$ satisfies $\ttau_l(\tL) = \tL[1-n]$ for
any $\LL^\infty$-grading $\tL$ of $L$.
\end{lemma}

\proof Clearly it is enough to prove the corresponding fact for the local model
$\ttau$. We find it convenient to use complex coordinates $U = \{ \xi \in \C^{n+1}
\suchthat |\re \xi| \leq 1, |\im\; \xi| = 1, \leftsc \re \xi,\im\; \xi \rightsc_{\R} =
0\}$. The tangent spaces are $TU_\xi = \{\eta \in \C^{n+1} \suchthat \leftsc
\im\;\xi,\im\;\eta \rightsc_{\R} = 0, \; \o_{\C^{n+1}}(\bar{\xi},\eta) = 0\}$. Hence if
$\Lambda$ is a Lagrangian subspace of $TU_\xi$, $\Lambda \oplus \R\bar{\xi}$ is a
Lagrangian subspace of $\C^{n+1}$. This stabilization defines a map $r: \LL
\longrightarrow \LL(2n+2)$ whose restriction to any fibre $\LL_\xi$ induces an
isomorphism of the fundamental groups. It follows that as our Maslov covering
$\LL_U^\infty$, we can take the pullback $r^*(\LL^\infty(2n+2))$ of the universal cover
of $\LL(2n+2)$.

$\tau^2$ is the time-one map of the Hamiltonian flow $\phi_t(\xi) =
\sigma_{t(2\psi(|\xi|)-1)}(\xi)$. Lift $(\phi_t)$ to an isotopy $(\tphi_t)$ of
$\LL_U^\infty$-graded symplectic automorphisms, starting with $\tphi_0 = \id$. By
definition, $\phi_t$ is equal to the identity in a neighbourhood of the zero-section
$S^n \subset U$ for any $t$. This implies that $\tphi_t(\tilde{S}^n) = \tilde{S}^n$ for
any $t$ and any $\LL_U^\infty$-grading $\tilde{S}^n$ of $S^n$. On the other hand,
$\tphi_1$ acts as some shift $[k]$ on the part of $\LL^\infty$ which lies over
$\partial U$. Using the fact that $\LL_U^\infty$ is defined as a pullback, it is not
difficult to see that $k$ can be computed as follows: choose a point $\xi \in \partial
U$ and a Lagrangian subspace $\Lambda \subset TU_\xi$. Then $\lambda(t) =
r(D\phi_t(\Lambda))$ is a loop in $\LL(2n+2)$, and one has
\[
k = -\leftsc C(2n+2),[\lambda] \rightsc.
\]
To determine $k$ it is convenient to take $\Lambda$ tangent to $\partial U$, because
$\phi_t|\partial U$ agrees with the inverse of the standard diagonal circle action on
$\C^{n+1}$. Then $\lambda(t) = e^{-2\pi i t}\Lambda \oplus e^{2\pi i t}(\R \bar{\xi})$,
which means that $k = 2n-2$. We now know that $\tphi_1 \circ [2-2n]$ is an
$\LL^\infty_U$-grading of $\tau^2$ which acts trivially on the boundary of $U$.
Therefore it must be the square $\ttau^2$ of the preferred grading $\ttau$. Hence
$\ttau^2(\tilde{S}^n) = \tphi_1(\tilde{S}^n[2-2n]) = \tilde{S}^n[2-2n]$ which proves
that $\ttau(\tilde{S}^n) = \tilde{S}^n[1-n]$. \qed

The next two results are graded analogues of Lemma \ref{th:pre-braiding} and Lemma
\ref{th:second-pre-braiding}.

\begin{lemma} \label{th:graded-pre-braiding}
Assume that $(l_1,l_2)$ form an $(A_2)$-configuration in $M$. Choose
$\LL^\infty$-gradings $\tL_1$, $\tL_2$ such that the index at the only intersection
point $\{x_0\} = L_1 \cap L_2$ is $\tilde{I}(\tL_1,\tL_2;x_0) = 1$. Then
$\ttau_{l_1}(\tL_2)$ and $\ttau^{-1}_{l_2}(\tL_1)$ are isotopic as $\LL^\infty$-graded
Lagrangian submanifolds.
\end{lemma}

\proof Let $\Sigma = L_1 \# L_2$ be the Lagrangian surgery, and $\tilde{\Sigma}$ the
$\LL^\infty$-grading from Lemma \ref{th:graded-surgery}. Recall that $\tau_{l_1}(L_2)$
and $\Sigma$ agree outside a neighbourhood of $L_1$. The gradings $\ttau_{l_1}(\tL_2)$
and $\tSigma$ also agree there: this follows from the fact that $\ttau_{l_1}$ is
trivial away from $L_1$, and that the grading $\tSigma$ agrees with $\tL_2$ on $\Sigma
\cap L_2$. In the proof of Lemma \ref{th:pre-braiding} we have used an isotopy from
$\tau_{l_1}(L_2)$ to $\Sigma$ which is concentrated near $L_1$. By considering a point
which remains fixed, one sees that when one lifts this to an isotopy of graded
Lagrangian submanifolds, it connects $\ttau_{l_1}(\tL_2)$ with $\tSigma$. The same kind
of argument shows that $\tSigma$ is isotopic to $\ttau_{l_2}^{-1}(\tL_1)$. \qed

\begin{lemma} \label{th:graded-second-pre-braiding}
Assume that $(l_1,l_2)$ is an $(A_2)$-configuration in $M$, and define $\tilde{g} =
(\ttau_{l_1}\ttau_{l_2})^3$. Then for any $\LL^\infty$-gradings $\tL_1$, $\tL_2$ one
has
\[
\tilde{g}(\tL_1) \grisotopic \tL_1[4-3n], \quad \tilde{g}(\tL_2) \grisotopic
\tL_2[4-3n],
\]
where $\grisotopic$ stands for isotopy of $\LL^\infty$-graded Lagrangian submanifolds.
\end{lemma}

\proof Clearly, if the result holds for some grading of $L_1,L_2$ then it holds for all
gradings. Hence we can assume that the absolute index $\tilde{I}(\tL_1,\tL_2;x_0)$ at
$\{x_0\} = L_1 \cap L_2$ is $1$. Lemma \ref{th:graded-pre-braiding} says that
$\ttau_{l_1}(\tL_2) \grisotopic \ttau_{l_2}^{-1}(\tL_1)$. Applying the same result to
$\tL_2$ and $\tL_1[n-2]$, which satisfy $\tI(\tL_2,\tL_1[n-2];x_0) = n -
\tI(\tL_1[n-2],\tL_2;x_0) = n - \tI(\tL_1,\tL_2;x_0) - n + 2 = 1$, yields
$\ttau_{l_2}(\tL_1) \grisotopic \ttau_{l_1}^{-1}(\tL_2)[2-n]$. Using these two
equations and Lemma \ref{th:preferred-grading} one computes
\begin{align*}
(\ttau_{l_1}\ttau_{l_2})^3(\tL_2) &= (\ttau_{l_1}\ttau_{l_2})^2\ttau_{l_1}(\tL_2)[1-n]
\grisotopic (\ttau_{l_1}\ttau_{l_2})^2\ttau_{l_2}^{-1}(\tL_1)[1-n] = \\ &=
\ttau_{l_1}\ttau_{l_2}(\tL_1)[2-2n] \grisotopic
\ttau_{l_1}\ttau_{l_1}^{-1}(\tL_2)[4-3n] = \tL_2[4-3n].
\end{align*}
The result for $\tL_1$ is proved in the same way. \qed

\subsection{Symplectically knotted Lagrangian spheres\label{subsec:knotted}}

After these preparations, we can now apply the `graded' techniques to the construction
of symplectically knotted Lagrangian spheres. The two parallel results obtained in this
way (for even-dimensional and odd-dimensional spheres, respectively) are

\begin{theorem} \label{th:knotted-even}
Let $(M^{2n},\o)$ be a compact symplectic manifold with contact type boundary, with $n$
even, which satisfies $[\o] = 0$ and $2c_1\mo = 0$. Assume that $M$ contains an
$(A_3)$-configuration $(l_0,l_1,l_2)$ of Lagrangian spheres. Then $M$ contains
infinitely many symplectically knotted Lagrangian spheres. More precisely, if one
defines $L_1^{(k)} = \tau_{l_2}^{2k}(L_1)$ for $k \in \Z$, then all the $L_1^{(k)}$ are
isotopic as smooth submanifolds of $M$, but no two of them are isotopic as Lagrangian
submanifolds.
\end{theorem}

\begin{theorem} \label{th:knotted-odd}
Let $(M^{2n},\o)$ be a compact symplectic manifold with contact type boundary, with $n
\geq 5$ odd, which satisfies $[\o] = 0$ and $2c_1\mo = 0$. Assume that $M$ contains an
$(A_4)$-configuration $(l_0,l_1,l_2,l_3)$ of Lagrangian spheres. Then $M$ contains
infinitely many symplectically knotted Lagrangian spheres. More precisely, if one
defines $L_1^{(k)} = g^{2k}(L_1)$ for $k \in \Z$, where $g = (\tau_{l_2}\tau_{l_3})^3$,
then all the $L_1^{(k)}$ are isotopic as smooth submanifolds of $M$, but no two of them
are isotopic as Lagrangian submanifolds.
\end{theorem}

Examples of manifolds satisfying these conditions are the Milnor fibres \eqref{eq:ak}.
Together these two Theorems prove the existence of infinitely many knotted Lagrangian
$n$-spheres in any dimension $n \neq 1,3$. As mentioned in the Introduction, the case
$n = 2$ has been proved before in \cite{seidel98b}, and the proof given there would
work in all even dimensions. Nevertheless, our new proof is substantially simpler.

\proof[Proof of Theorem \ref{th:knotted-even}] Let $\{x_0\} = L_0 \cap L_1$ and
$\{x_1\} = L_1 \cap L_2$. Choose an $\infty$-fold Maslov covering $\LL^\infty$ on $M$,
and $\LL^\infty$-gradings $\tL_0,\tL_1,\tL_2$ in such a way that $\tI(\tL_0,\tL_1;x_0)
= \tI(\tL_1,\tL_2;x_1) = 0$. Set $\tL_1^{(k)} = \ttau_{l_2}^{2k}(\tL_1)$. The
assumptions on $M$ imply that Floer cohomology (and its graded version) are
well-defined for all (graded) Lagrangian submanifolds of $M$ with vanishing first Betti
number. Because $L_0 \cap L_1$ and $L_1 \cap L_2$ intersect in a single point, which
has index zero, it follows from the definition of graded Floer cohomology that
\[
HF^*(\tL_0,\tL_1) = \Z/2^{[0]}, \quad HF^*(\tL_1,\tL_2) = \Z/2^{[0]}.
\]
Here $\Z/2^{[k]}$ stands for the $\Z$-graded group which is $\Z/2$ in degree $k$ and
zero in other degrees. Since $L_0 \cap L_2 = \emptyset$, one can choose $\tau_{l_2}$ in
such a way that it acts trivially in a neighbourhood of $L_0$. This implies that $L_0
\cap L_1^{(k)} = L_0 \cap L_1 = \{x_0\}$. Because the preferred grading $\ttau_{l_2}$
also acts trivially near $L_0$, one has $\tI(\tL_0,\tL_1^{(k)};x_0) =
\tI(\tL_0,\tL_1;x_0) = 0$ and hence (for any $k$)
\begin{equation} \label{eq:first-reference}
HF^*(\tL_0,\tL_1^{(k)}) = \Z/2^{[0]}.
\end{equation}
On the other hand, using the invariance of graded Floer cohomology under graded
symplectic automorphisms together with Lemma \ref{th:preferred-grading}, one finds that
\begin{equation} \label{eq:second-reference}
\begin{aligned}
HF^*(\tL_1^{(k)},\tL_2) & \iso HF^*(\tL_1,\ttau_{l_2}^{-2k}(\tL_2)) =\\&=
HF^*(\tL_1,\tL_2[2k(n-1)]) = \Z/2^{[2k(1-n)]}. \end{aligned}
\end{equation}
Now assume that for some $L_1^{(k)}$, for some $k$, is Lagrangian isotopic to $L_1$.
This implies that $\tL_1^{(k)}$ is graded Lagrangian isotopic to $\tL_1[r]$ for some $r
\in \Z$. Using the isotopy invariance property of Floer cohomology one obtains
\[
HF^*(\tL_0,\tL_1^{(k)}) \iso HF^*(\tL_0,\tL_1[r]) = \Z/2^{[-r]} \text{ and }
HF^*(\tL_1^{(k)},\tL_2) \iso \Z/2^{[r]}.
\]
Comparing the first part of this equation with \eqref{eq:first-reference} yields $r =
0$. Comparing the second part with \eqref{eq:second-reference} yields $r = 2k(1-n) = 0$
and hence $k = 0$. This means that $L_1^{(k)}$ is not Lagrangian isotopic to $L_1$
unless $k = 0$. Because of the way in which the $L_1^{(k)}$ are defined, it follows
that no two of them are Lagrangian isotopic. The topological part of the theorem
follows from Lemma \ref{th:haefliger}\ref{item:even}. \qed

\proof[Proof of Theorem \ref{th:knotted-odd}] Since this is very similar to the proof
of Theorem \ref{th:knotted-even}, we will be more brief. Take an $\infty$-fold Maslov
covering $\LL^\infty$ and gradings $\tL_0,\tL_1,\tL_2$ such that $HF^*(\tL_0,\tL_1) =
HF^*(\tL_1,\tL_2) = \Z/2^{[0]}$. Define $\tilde{g} = (\ttau_{l_2}\ttau_{l_3})^3$ and
$\tL_1^{(k)} = \tilde{g}^{2k}(\tL_1)$. One can choose $\tau_{l_2},\tau_{l_3}$ in such a
way that $g$ and $\tilde{g}$ act trivially near $L_0$. This implies that
$HF^*(\tL_0,\tL_1^{(k)}) = \Z/2^{[0]}$ for any $k$. Using Lemma
\ref{th:graded-second-pre-braiding} one computes
$
HF^*(\tL_1^{(k)},\tL_2) \iso HF^*(\tL_1,\tilde{g}^{-2k}(\tL_2)) \iso
HF^*(\tL_1,\tL_2[2k(3n-4)]) = \Z/2^{[2k(4-3n)]}$. The assumption that $\tL_1^{(k)}$,
for some $k \neq 0$, is graded Lagrangian isotopic to $\tL_1[r]$, for some $r \in \Z$,
leads to the contradiction $0 = r = 2k(4-3n)$. The topological part of the theorem is
Lemma \ref{th:haefliger}\ref{item:odd}. \qed

\subsection{$K3$ and Enriques surfaces\label{subsec:k3}}

Let $\mo$ be a closed symplectic four-manifold such that $c_1(M)$ is a torsion class.
We want to consider Lagrangian spheres in $M$. This is a borderline case for Floer
cohomology, where the conventional methods do not yield a completely satisfactory
theory. The problems appear when one tries to prove that the Floer group is independent
of the choice of almost complex structure. We will now explain what parts of Floer's
construction can be salvaged from this breakdown. Throughout the whole of this section,
{\em all Lagrangian submanifolds are assumed to be two-spheres}. For such an $L \subset
M$, let $\JJ^{reg}(L)$ be the set of $\o$-compatible almost complex structures $J$ such
that there are no non-constant $J$-holomorphic maps $\CP{1} \rightarrow M$ or
$(D^2,\partial D^2) \rightarrow (M,L)$.

\begin{lemma} $\JJ^{reg}(L)$ is a dense subset of the space of all $\o$-compatible almost complex
structures.
\end{lemma}

\proof The virtual dimension of the space of $J$-holomorphic spheres in a homology
class $A$ is $4+2\leftsc c_1(M), A \rightsc - 6 = -2$. Similarly, the virtual dimension
of the space of $J$-holomorphic discs representing $B \in H_2(M,L)$ is $2+\leftsc
2c_1(M,L), B \rightsc - 3 = -1$. Here $2c_1(M,L) \in H^2(M,L)$ is the relative first
Chern class, and we have used the assumption that $H^1(L) = 0$ to conclude that
$c_1(M,L)$ is a torsion class. Standard transversality results say that for generic $J$
there are no {\em simple} $J$-holomorphic spheres and no {\em somewhere injective}
$J$-holomorphic discs. The first result implies that there are no non-constant
$J$-spheres at all, because any such sphere covers a simple one. Similarly, the second
result implies that there are no $J$-holomorphic discs. However, this time the argument
is more complicated: it relies on the structure theorem of Kwon-Oh \cite{oh-kwon96} or
on the simpler form given by Lazzarini \cite{lazzarini98}. \qed

Let $L_0,L_1 \subset M$ be two Lagrangian submanifolds which intersect transversally.
Fix $J_0 \in \JJ^{reg}(L_0)$, $J_1 \in \JJ^{reg}(L_1)$. We define
$\JJ^{reg}(L_0,J_0,L_1,J_1)$ to be the space of smooth families $\J = (J_t)_{0 \leq t
\leq 1}$ of compatible almost complex structures, connecting the given $J_0,J_1$, which
satisfy the following conditions: (a) there are no non-constant $J_t$-holomorphic maps
$\CP{1} \rightarrow M$ for any $t$; (b) any solution of Floer's equation, $u \in
\moduli(x_-,x_+;\J)$ for $x_-,x_+ \in L_0 \cap L_1$, is regular.

\begin{lemma}
$\JJ^{reg}(L_0,J_0,L_1,J_1)$ is a dense subset of of the space of all families $\J$
which connect $J_0$ with $J_1$.
\end{lemma}

The proof is again a combination of standard transversality arguments. For $\J \in
\JJ^{reg}(L_0,J_0,L_1,J_1)$ one can define the Floer cohomology
$HF(L_0,J_0,L_1,J_1;\J)$ in the familiar way, using a suitable Novikov ring $\Lambda$
as coefficient ring. The next step in the construction is

\begin{lemma} \label{th:half-independent}
$HF(L_0,J_0,L_1,J_1;\J)$ is independent of $\J \in \JJ^{reg}(L_0,J_0,L_1,J_1)$ up to
canonical isomorphism.
\end{lemma}

The proof uses the continuation equation $\partial_su + \hat{J}_{s,t}(u)\partial_tu =
0$ for a two-parameter family $\hat{\J} = (\hat{J}_{s,t})$ of almost complex structures
such that $\hat{J}_{s,0} = J_0$ and $\hat{J}_{s,1} = J_1$ for all $s \in \R$. There
will be $\hat{J}_{s,t}$-holomorphic spheres (for isolated values of $s,t$) even if
$\hat{\J}$ is generic, but these can be dealt with as in \cite{hofer-salamon95}.
Following the usual custom, we will from now on omit $\J$ from the notation of Floer
homology.

The problematic issue mentioned above is whether $HF(J_0,L_0,J_1,L_1)$ is independent
of $J_0,J_1$. We will use only a special case, in which the result is obvious:

\begin{lemma} \label{th:one-point}
Assume that $L_0,L_1 \subset M$ intersect transversally in a single point. Then
$HF(L_0,J_0,L_1,J_1) \iso \Lambda$ for all $J_0 \in \JJ^{reg}(L_0)$, $J_1 \in
\JJ^{reg}(L_1)$. \qed
\end{lemma}

By definition, Floer cohomology is invariant under symplectic automorphisms, in the
sense that $HF(\phi(L_0),\phi_*(J_0),\phi(L_1),\phi_*(J_1)) \iso HF(L_0,J_0,L_1,J_1)$.
The next result is a weak form of isotopy invariance.

\begin{lemma} \label{th:pre-isotopy-invariance}
Assume that $L_0,L_1$ are transverse and choose $J_0 \in \J^{reg}(L_0)$, $J_1 \in
\J^{reg}(L_1)$. Let $\phi \in \Aut\mo$ be a map which is Hamiltonian isotopic to the
identity, and such that $\phi(L_1)$ intersects $L_0$ transversally. Then $\phi_*(J_1)
\in \JJ^{reg}(\phi(L_1))$ and
\[
HF(L_0,J_0,\phi(L_1),\phi_*(J_1)) \iso HF(L_0,J_0,L_1,J_1).
\]
\end{lemma}

\proof[Outline of the proof.] take a function $H \in \smooth([0;1] \times M,\R)$ with
$H_t = 0$ for $t \leq 1/3$ or $t \geq 2/3$, such that the Hamiltonian isotopy
$(\phi_t^H)$ generated by it satisfies $\phi_1^H = \phi^{-1}$. Let $X_t$ be the
Hamiltonian vector field of $H_t$, and $\psi \in \smooth(\R,\R)$ a cutoff function with
$\psi(s) = 0$ for $s \leq 0$ and $\psi(s) = 1$ for $s \geq 1$. One considers finite
energy solutions $u: \R \times [0;1] \longrightarrow M$ of the equation
\begin{equation} \label{eq:continuation-two}
\begin{cases} & \partial_su + \hat{J}_{s,t}(u)(\partial_tu - \psi(s)X_t(u)) = 0, \\
& u(s,0) \in L_0, \quad u(s,1) \in L_1. \end{cases}
\end{equation}
Here $\hat{\J} = (\hat{J}_{s,t})$ is a two-parameter family of $\o$-compatible almost
complex structures such that $\hat{J}_{s,0} = J_0$, $\hat{J}_{s,1} = J_1$ for all $s$.
If one writes $u(s,t) = \phi_t^H(w(s,t))$ then \eqref{eq:continuation-two} in the
region $s \geq 1$ reads
\[
\begin{cases} & \partial_sw + \hat{J}'_{s,t}(w)\partial_tw = 0, \\
& w(s,0) \in L_0, \quad w(s,1) \in \phi(L_1). \end{cases}
\]
where $\hat{J}_{s,t}' = (\phi_t^H)^{-1}_*\hat{J}_{s,t}$, in particular $\hat{J}_{s,1}'
= \phi_*J_1$. Following the usual strategy, one can use solutions of
\eqref{eq:continuation-two} to define a map
\begin{equation} \label{eq:continuation-three}
HF(L_0,J_0,L_1,J_1) \longrightarrow HF(L_0,J_0,\phi(L_1),\phi_*(J_1)).
\end{equation}
The main technical point is that there can be no bubbling off of holomorphic discs,
since the almost complex structures $\hat{J}_{s,t}$ for $t = 0,1$ do not admit such
discs. Standard arguments of a similar kind show that \eqref{eq:continuation-three} is
an isomorphism. \qed

If $2c_1(M) = 0$, one can consider graded Lagrangian submanifolds, and obtains graded
Floer groups $HF^*(\tL_0,J_0,\tL_1,J_1)$ with properties analogous to those above.

\begin{theorem} \label{th:k3}
Let $\mo$ be a symplectic four-manifold with $2c_1\mo = 0$ and which contains an
$(A_3)$-configuration $(l_0,l_1,l_2)$ of Lagrangian two-spheres. Define $L_1^{(k)} =
\tau_{l_2}^{2k}(L_1)$. Then all the $L_1^{(k)}$ are isotopic as smooth submanifolds,
but no two of them are isotopic as Lagrangian submanifolds.
\end{theorem}

\proof This is the same argument as in Theorem \ref{th:knotted-even}, except that one
has to be more careful about the properties of Floer cohomology. Choose an
$\infty$-fold Maslov covering $\LL^\infty$ and gradings $\tL_k$ such that
$\tI(\tL_0,\tL_1;x_0) = \tI(\tL_1,\tL_2;x_1) = 0$, where $x_0,x_1$ are the unique
intersection points. Fix some $k \neq 0$ and write $L_1' = \tau_{l_2}^{2k}(L_1)$,
$\tL_1' = \ttau_{l_2}^{2k}(\tL_1)$. Assume that $L_1 \isotopic L_1'$, so that $\tL_1[r]
\grisotopic \tL_1'$ for some $r \in \Z$. One can embed the Lagrangian isotopy into a
Hamiltonian isotopy $(\phi_t)$ of $M$. The graded analogue of Lemma
\ref{th:pre-isotopy-invariance} shows that
\begin{align*}
HF^*(\tL_0,J_0,\tL_1,J_1) &\iso HF^{*-r}(\tL_0,J_0,\tL_1',(\phi_1)_*(J_1)), \\
HF^*(\tL_1,J_1,\tL_2,J_2) &\iso HF^{*+r}(\tL_1',(\phi_1)_*(J_1),\tL_2,J_2)
\end{align*}
for all $J_m \in \JJ^{reg}(L_m)$, $m = 0,1,2$. On the other hand, by using the graded
version of Lemma \ref{th:one-point} and arguing as in the proof of Theorem
\ref{th:knotted-even} one finds that
\begin{align*}
& HF^*(\tL_0,J_0,\tL_1,J_1) = HF^*(\tL_1,J_1,\tL_2,J_2) = \Lambda^{[0]}, \\ &
HF^*(\tL_0,J_0,\tL_1',(\phi_1)_*(J_1)) = \Lambda^{[0]}, \quad
HF^*(\tL_1',(\phi_1)_*(J_1),\tL_2,J_2) = \Lambda^{[-2k]},
\end{align*}
which leads to a contradiction. \qed

\begin{cor}
For $\mo$ as in Theorem \ref{th:k3}, the map $\pi_0(\Aut\mo) \longrightarrow
\pi_0(\Diff(M))$ has infinite kernel.
\end{cor}

\proof We have already quoted the fact that $\tau_{l_2}^2$ is trivial in
$\pi_0(\Diff(M))$ \cite[Lemma 6.3]{seidel98b}. Theorem \ref{th:k3} obviously implies
that $[\tau_{l_2}^2] \in \pi_0(\Aut\mo)$ has infinite order. \qed

\begin{example} \label{ex:enriques}
Recall that an Enriques surface is an algebraic surface with fundamental group $\Z/2$,
whose universal cover is a $K3$ surface. Enriques surfaces satisfy $2c_1 = 0$. We will
now construct an Enriques surface which contains an $(A_3)$-configuration of Lagrangian
spheres. Consider the quartic surface $X \subset \CP{3}$ defined by
\begin{equation} \label{eq:singular-quartic}
(x_0^2 + x_1^2)^2 + x_0^2x_2^2 + x_2^4 + x_3^4 = 0.
\end{equation}
$X$ has two singular points $[1:\pm i:0:0]$ of type $(A_3)$. Let $Y$ be the minimal
resolution of singularities of $X$. It is a $K3$ surface and hence has a holomorphic
symplectic form $\O$. Complex conjugation defines an anti-holomorphic involution on
$X$, which is free because \eqref{eq:singular-quartic} has no nonzero real solutions.
Because of the uniqueness of minimal resolutions, this involution lifts to an
anti-holomorphic involution on $Y$, which we denote by $\iota$. Obviously $\iota$ is
again free. Since the holomorphic symplectic form is unique up to a constant, it
satisfies $\iota^*\O = z\bar{\O}$ for some $z \in S^1$. By rescaling $\Omega$ we can
assume that $z = 1$. Then $\o = \re \Omega$ descends to a real symplectic form on the
quotient $M = Y/\iota$. The singular points of $X$ give rise to two disjoint
$(A_3)$-configuration of rational curves in $Y$, which are exchanged by $\iota$. These
curves are Lagrangian with respect to $\o$, hence descend to a single
$(A_3)$-configuration of Lagrangian spheres on $M$. To see that $M$ is K{\"a}hler one
argues as follows: let $\beta \in \Omega^{1,1}$ be a positive form on $Y$ such that
$\iota^*\beta = -\beta$. Such a form can be constructed e.g. by averaging. Let $g$ be
the unique Ricci-flat K{\"a}hler metric which has the same K{\"a}hler class as $\beta$.
The uniqueness theorem for such metrics implies that $g$ is $\iota$-invariant. The
metric $g$ is hyperk{\"a}hler, which means that there are complex structures $J$ and
$K$ such that $\O(v,w) = g(v,Jw) + ig(v,Kw)$. It follows that $J$ is $\iota$-invariant
and hence descends to a complex structure on $M$ which is compatible with $\o$. This
means that $M$ is K{\"a}hler and in fact an Enriques surface.
\end{example}

\begin{example} \label{ex:k3}
$K3$ surfaces have $c_1 = 0$. As an example of a $K3$ surface containing an
$(A_3)$-configuration of Lagrangian spheres one can take the manifold $(Y,\re \Omega)$
constructed in the previous example.
\end{example}

\providecommand{\bysame}{\leavevmode\hbox to3em{\hrulefill}\thinspace}

\end{document}